\documentclass[a4paper]{amsart}
\usepackage{amsfonts}
\usepackage{amsmath}
\usepackage{mathrsfs}
\usepackage{graphicx}
\usepackage{amsmath,amsthm,amssymb,amscd}
\usepackage{color}
\usepackage{enumerate}

\usepackage{epsfig}

\theoremstyle{plain}\newtheorem{maintheorem}{Theorem}

\newtheorem{maincorollary}{Corollary}

\newtheorem{mainlemma}{Lemma}

\newtheorem{Thm}{Theorem}[section]
\newtheorem{Lem}[Thm]{Lemma}
\newtheorem{Prop}[Thm]{Proposition}

\theoremstyle{remark}
\newtheorem{Def}[Thm] {Definition}
\newtheorem{Rem}[Thm] {Remark}

\newtheorem{Que}[Thm] {Question}
\newtheorem{thm}{Theorem}[section]
\newtheorem{lem}[thm]{Lemma}

\theoremstyle{remark}

\theoremstyle{definition}

\newcommand{\htop}{h_{\text{top}}}

\long\def\begcom#1\endcom{}

\newcommand{\length}{\operatorname{\length}}

\def\length{\operatorname{length}}

\newcommand{\bl} {\begin{lemma}}
\newcommand{\el} {\end{lemma}}
\newcommand{\bt} {\begin{theorem}}
\newcommand{\et} {\end{theorem}}

\newcommand{\bp}{\begin{proof}}
\newcommand{\ep}{\end{proof}}
\newcommand {\be}{\begin{equation}}
\newcommand  {\ee} {\end{equation}}
\newcommand  {\beq} {\begin{eqnarray*}}
\newcommand  {\eeq} {\end{eqnarray*}}

\newcommand  {\bd} {\begin{definition}}
\newcommand  {\ed} {\end{definition}}

\newtheorem{theorem}{Theorem}[section]
\newtheorem{lemma}[theorem]{Lemma}

\theoremstyle{definition}
\newtheorem{definition}[theorem]{Definition}

\theoremstyle{remark}

\numberwithin{equation}{section}

\begin{document}


\title[Transitively-Saturated Property, Banach Recurrence and Lyapunov Regularity ] {Transitively-Saturated Property, Banach Recurrence and Lyapunov Regularity }
\author[Y. Huang] {Yu Huang}
\address[Y. Huang]{Department of Mathematics, Sun Yat-sen University, Guangzhou 510275, People¡¯s Republic of China}
\email{stshyu@mail.sysu.edu.cn }

\author[X. Tian] {Xueting Tian}
\address[X. Tian]{School of Mathematical Sciences,  Fudan University   \\  Shanghai 200433, People's Republic of China}
\email{xuetingtian@fudan.edu.cn}

\author[X. Wang] {Xiaoyi Wang}
\address[X. Wang]{College of Mathematics and Statistics, Zhaoqing University,
 Zhaoqing 526061, GuanDong,
 People's Republic of China  }
\email{wxy635@163.com}

\maketitle

\def\abstractname{\textbf{Abstract}}

\begin{abstract}\addcontentsline{toc}{section}{\bf{English Abstract}}

The topological entropy of various  gap-sets on periodic-like recurrence and Birkhoff regularity were considered in \cite{T16} but some Banach recurrence and Lyapunov regularity are not considered.
 In this paper we introduce five new levels  on Banach recurrence and show they all carry full topological entropy,   and   simultaneously combine  with Lyapunov regularity to get some refined theory on mixed  multifractal analysis of \cite{Barreira-S2001,FengH}.


   In this process, we strengthen Pfister and Sullivan's result of \cite{PS} from saturated property to transitively-saturated property (and from single-saturated property to transitively-convex-saturated property).

\end{abstract}

\tableofcontents

\section{Introduction} \setlength{\parindent}{2em}
Many periodic-like recurrence such as periodic, almost periodic, weakly almost periodic, quasi-weakly almost periodic and Birkhoff regularity were considered in \cite{T16} in the viewpoint of
 topological  entropy. But some Banach recurrence and Lyapunov regularity are not considered.
 In this paper we introduce five new levels  on Banach recurrence and show they all carry full topological entropy,
  and   simultaneously combine  with Lyapunov regularity to get some refined theory on mixed  multifractal analysis
   of \cite{Barreira-S2001,FengH}. We refer to \cite{Walter} for ergodic theory and to \cite{Bar,Pesin97-book} for dimension theory. 
 In this paper, a dynamical system $(X,  T )$ means always that $X$ is a compact metric
space and $T:X\rightarrow X$ is a continuous map.  Let $M(T,X)$ denote the space of all $T-$invariant measures and  $M_{erg}(T,X)$ denote the space of all $T-$ergodic measures.

\subsection{Main Results}

\subsubsection{Transitively-saturated property}



Given $x\in X$, let $\omega_T(x)$ denote the $\omega-$limit set of $x$,    $M_x$ be the limit set of the empirical measures for $x$. Define
  $ Rec =\{x\in X\,|x\in \omega_T(x)  \}$ and  $Tran   = \{x\in X|   \omega_T(x) = X\}.$
 For $A\subseteq X,$ let $h_{top}(A)$ denote the topological entropy of $A$ defined by Bowen in  \cite{Bowen1} ( see Section \ref{section-prepare})
  and given an invariant measure $\mu,$ let  $h_\mu(T)$ denote its metric entropy of $\mu$ and let $S_\mu$ denote  the support of $\mu$..

 \begin{Def}\label{Def-saturated}
  We say that the system $T$ has {\it  transitively-saturated} property or $T$ is {\it  transitively-saturated}, if  for any  compact connected nonempty set $K \subseteq M (T,X ),$
\begin{eqnarray} \label{eq- saturated-definition}
 G^T_K\neq \emptyset \,\,\text{ and }\,\,h_{top} (T,G^T_K )=\inf\{h_\mu (T)\,|\,\mu\in K\},
\end{eqnarray}  where $G^T_{K}
=\{x\in Tran|\, M_x=K
  \}.$ 
     If above equality (\ref{eq- saturated-definition}) only holds for the $K$ that is contained in some   convex sum of finite  invariant measures, we say the system $T$    is    transitively-convex-saturated. If above equality (\ref{eq- saturated-definition}) only holds for the $K$ that consists of one invariant measure, we say the system $T$    is    transitively-single-saturated.

    \end{Def}

    \begin{maintheorem}\label{Thm--PropMain-productproperty-uniformseparation}
Suppose that $(X,\,T)$ satisfies $g$-almost product property and there is an invariant measure
with full support.
 Then \\
 (1)
   $T$ is  transitively-convex-saturated, the set  $$\{\omega\in M(T,X)|\,\mu\text{ is ergodic, } \,S_\omega\text{ is minimal }\}$$
  is dense in $M(T,X)$  and almost periodic set $A$ is dense in $X$. \\
  (2)
    If further $T$ has uniform separation, then $T$ is  transitively-saturated.

 \end{maintheorem}

Theorem \ref{Thm--PropMain-productproperty-uniformseparation} generalizes Pfister and Sullivan's result of \cite{PS} from saturated property (for which they used $G_K$ to replace $G^T_K $     in  (\ref{eq- saturated-definition})) to transitively-saturated property. We will give the proof of Theorem \ref{Thm--PropMain-productproperty-uniformseparation} in Section \ref{Section-Proposition-MAIN}.
\begin{Rem}
By    \cite[Theorem 3]{Bowen1}  if $K$ is a singleton only with one ergodic measure $\mu,$ then one has naturally \begin{eqnarray}
G_K\neq \emptyset,\,\,\mu(G_K)=1, \,\,\text{ and }\,\, h_{top} (T,G_K)= h_\mu(T).
\end{eqnarray}  On the other hand, by \cite[Theorem 3]{Bowen1} $h_{top} (T,\Gamma)\geq h_\mu(f)$ for any $\Gamma$ with $\mu(\Gamma)=1$ so that $h_{top} (T,G_K\cap S_\mu)= h_\mu(f)$. This holds for    any dynamical system $(X,T)$ and any ergodic measure $\mu.$  If $S_\mu=X,$ then $G_K=G_K^T$ so that $h_{top} (T,G_K^T)=h_{top} (T,G_K) = h_\mu(T)$. However, in general
 it is not true that  $h_{top} (T,G^T_K)=h_{top} (T,G_K\setminus S_\mu)= h_\mu(f)$ when $S_\mu\neq X$.
  For example, if $T$ is not transitive with positive topological entropy (e.g., a product map by identity map on $[0,1]$ and a system with positive topological entropy), then for any $K=\{\mu\}$ that $\mu$ is ergodic with positive metric entropy, $G^T_K$ is empty so that  $h_{top} (T,G^T_K)=0< h_\mu(f).$
\end{Rem}

\subsubsection{Banach Recurrence and full entropy gap}

Let $S\subseteq \mathbb{N}$,   define
$$\bar{d} (S):=\limsup_{n\rightarrow\infty}\frac{|S\cap \{0,  1,  \cdots,  n-1\}|}n,  \,  \,   \underline{d} (S):=\liminf_{n\rightarrow\infty}\frac{|S\cap \{0,  1,  \cdots,  n-1\}|}n,  $$
where $|Y|$ denotes the cardinality of the set $Y$.   These two concepts are called {\it upper density} and {\it lower density} of $S$,   respectively.   If $\bar{d} (S)=\underline{d}(S)=d,  $ we call $S$ to have density of $d.  $  Define $$B^* (S):=\limsup_{|I|\rightarrow\infty}\frac{|S\cap I|}{|I|},  \,  \,   B_* (S):=\liminf_{|I|\rightarrow\infty}\frac{|S\cap I|}{|I|},  $$
here $I\subseteq \mathbb{N}$ is taken from finite continuous integer intervals.   These two concepts are called {\it Banach upper density} and {\it Banach lower density} of $S$,   respectively.
    These concepts of density are basic and have played important roles in the field of dynamical systems,   ergodic theory and number theory,   etc.
Let $U,  V\subseteq X$ be two nonempty open subsets and $x\in X.  $
 Define
  sets of visiting time
 $$N (U,  V):=\{n\geq 1|\,  U\cap f^{-n} (V)\neq \emptyset\} \,  \,  \text{ and } \,  \,  N (x,  U):=\{n\ge 1|\,  f^n (x)\in U\}.  $$

\begin{Def}
 A point $x\in X$ is called Banach upper recurrent, if $N (x,  V_\epsilon (x)) $ has positive Banach upper density
 where $V_\epsilon (x)$ denotes  the  ball centered at $x$ with radius $\epsilon$. A
 point $x\in X$ is called     upper  recurrent (or quasi-weakly almost recurrent, see \cite{Zhou93}), if $N (x,  V_\epsilon (x)) $ has positive
 upper density.

 \end{Def}

Let $BR$ denote the set of all Banach upper recurrent points and let $QW$ denote
 the set of  upper recurrent points. The entropy estimate on classification of
 $QW$ has been studied in \cite{T16}.
  In present paper we will classify $BR\setminus QW$   and consider their entropy gap.  Given $x\in X,$   let $C_x=\overline{\cup_{m\in M_x  }S_m}$. Let $ BR^\# :=BR\setminus QW$,
 \begin{eqnarray*}   W^\# &:=&\{\, x\in  BR^\# \,|\,S_\mu=C_x  \text{ for every }\mu \in M_x \},
   \\   V^\#  &:=&\{ \,x\in  BR^\#  |\,\exists \,\mu\in M_x  \text{ such that } S_\mu=C_x\},\\
   S^\#  &:=&\{ \,x\in  X  |\,\cap_{\mu\in M_x} S_\mu\neq \emptyset\}.
          \end{eqnarray*}
           More precisely,  in present paper we are mainly to consider
 $ BR^\#$  which is divided into following several levels with
  different asymptotic behavior:
\begin{eqnarray*}   
   BR_1  &:=&  W^\#,\\
   BR_2  &:=& V^\#\cap S^\# ,\,\,\,\, BR_3   :=  V^\#, \\
   BR_4  &:=& V^\#\cup (BR^\#\cap S^\#) ,\,\,\,\,\, BR_5   :=  BR^\#,
          \end{eqnarray*} Note that $BR_1\subseteq BR_2\subseteq BR_3\subseteq BR_4\subseteq BR_5.$

\begin{Def}\label{def-full-entropygaps} For a collection of subsets $Z_1,Z_2,\cdots,Z_k\subseteq X$ ($k\geq 2$), we say $\{Z_i\}$ has {\it full entropy gaps} with respect to $Y\subseteq X$ (simply, FEG w.r.t $Y$) if
$$ h_{top} (T, (Z_{i+1}\setminus Z_i)\cap Y )=h_{top}(T, Y) \,\,\,\text{ for all } 1\leq i<k,$$
where $h_{top}(T,Z)$ denotes the topological entropy of a set $Z\subseteq X.$
 \end{Def}
Often, but not always, the sets $Z_i$ are nested ($Z_i\subseteq Z_{i+1}$).
 Note that for any system with zero topological entropy, it is obvious that any collection $\{Z_i\}$ has full entropy gaps with respect to any  $Y\subseteq X$.

\begin{Def}

We say $T$ satisfies the  entropy-dense property 
  if for any $\mu\in M(T,  X)$,   for any neighborhood $G$ of $\mu$ in $M(X)$,
   and for any $\eta>0$,   there exists a closed $T$-invariant set $\Lambda_{\mu}\subseteq X $ such that $M(T,  \Lambda_{\mu})\subseteq G$ and $\htop(T,  \Lambda_{\mu})>h_{\mu}-\eta$. By classical variational principle,
   it is equivalent that for any neighborhood $G$ of $\mu$ in $M(X)$,   and for any $\eta>0$,   there exists a $\nu\in M_{erg}(T,  X)$ such that $h_{\nu}>h_{\mu}-\eta$ and $M(T,  S_{\nu})\subseteq G$.


\end{Def}

 From \cite[Proposition 2.3 (1)]{PS2005}, entropy-dense property holds
  for dynamical systems with $g$-product property.
A point  $x\in X$ is generic for some invariant measure  $\mu$ means that $M_x =\{\mu\} $ (or equivalently,  Birkhoff averages of all
   continuous functions converge to the integral of $\mu$).    Let $G_\mu$   denote the set of all generic points for  $\mu $.
    Let $QR  = \cup_{\mu\in M(T,X)} G_\mu.$  The points in $QR$ are called quasiregular points of   $T$ in \cite{DGS,Oxt}.

\begin{maintheorem}\label{Maintheorem000-2017-banachrecurrence}
Suppose that $(X,\,T)$   has transitively-convex-saturated property and  entropy-dense property.
 If $(X,\,T)$ is not uniquely ergodic and  there is an invariant measure
with full support,    then
    $\{\emptyset, QR\cap BR_1 ,BR_1, BR_2,BR_3, BR_4,BR_5 \}$ has   full entropy gap w.r.t. $  Tran $.

 \end{maintheorem}

\subsubsection{Multifractal Analysis }

Pesin and Pitskel \cite{Pesin-Pitskel1984} are the first to notice the phenomenon of the irregular set carrying full topological entropy in the case of the full shift on two symbols from the dimensional perspective. Barreira, Schmeling, etc. studied   the irregular set  in the setting
of shifts of finite type and beyond, see \cite{Barreira-Schmeling2000,Bar,TV,Tho2012,To2010} etc. Ruelle uses
 the terminology in \cite{Ruelle} `historic behavior' to describe irregular point and in contrast to dimensional perspective,
  Takens asks in \cite{Takens} for which smooth dynamical systems the points with historic behavior has positive Lebesgue measure.
  One aim of present paper is to generalize these results into more general cases.
 Let $\alpha:M(T,X)\rightarrow \mathbb{R}$ be a continuous function, define
 $$I_\alpha:=\{x\in X: \inf_{\mu\in M_x}\alpha(\mu)<\sup_{\mu\in M_x}\alpha(\mu)\};$$
 $$R_\alpha:=\{x\in X: \inf_{\mu\in M_x}\alpha(\mu)=\sup_{\mu\in M_x}\alpha(\mu)\};$$
 $$R_\alpha(a):=\{x\in X: \inf_{\mu\in M_x}\alpha(\mu)=\sup_{\mu\in M_x}\alpha(\mu)=a\}.$$
  Let  $L_\alpha=[\inf_{\mu\in M_x}\alpha(\mu),\sup_{\mu\in M_x}\alpha(\mu)]$  and $Int(L_\alpha)$ denote its interior interval.   We list three conditions for $\alpha$:

 $[A.1].$ For any $\mu,\nu\in M(T,X)$,    $ \beta(\theta):=\alpha(\theta\mu+(1-\theta)\nu)$ is strictly monotonic on $[0,1]$ when $\alpha(\mu)\neq \alpha(\nu)$.

 $[A.2].$ For any $\mu,\nu\in M(T,X)$,   $ \beta(\theta):=\alpha(\theta\mu+(1-\theta)\nu)$ is constant on $[0,1]$ when $\alpha(\mu)= \alpha(\nu) $.

$[A.3].$ For any $\mu,\nu\in M(T,X)$,    $ \beta(\theta):=\alpha(\theta\mu+(1-\theta)\nu)$ is not constant over any subinterval of $[0,1]$ when $\alpha(\mu)\neq \alpha(\nu)$.

Note that $[A.1]$ implies $[A.3]$.

    \begin{Thm}\label{Thm--RoughResult-multifrcalanalysis}
 Let $\alpha:M(T,X)\rightarrow \mathbb{R}$ is a continuous function. 

 (1) Suppose that $(X,\,T)$   is  transitively-convex-saturated and $\alpha$ satisfies $[A.3]$, then either $I_\alpha $ is empty or
 $$h_{top}(T,I_\alpha\cap Tran)=h_{top}(T,I_\alpha )=h_{top}(T).$$


 (2) Suppose that $(X,\,T)$   is  transitively-single-saturated, 
  then for any $a\in\mathbb{R},$ either $R_\alpha(a)$ is empty or
  $$ h_{top} (T,R_\alpha(a)\cap Tran  )=h_{top} (T,R_\alpha(a)   )=\sup \{h_{\rho}  |\,\,\rho\in M (T,X)\,\,and\,\,\alpha (\rho)=a\} .$$

 \end{Thm}

 \begin{Rem}\label{Rem-Barreira-Saussol-mixed} Let $\phi,\psi$ be two continuous functions on $X$ and $\psi$ is required to be positive. Define $\alpha(\mu)=\frac{\int \phi d\mu}{\int \psi d\mu}$ which satisfies $[A.1]$, $[A.2]$ and $[A.3]$. In this case recall that  the part for $R_\alpha(a)$ has been discussed  in \cite{Barreira-S2001} for hyperbolic systems by Barreira and Saussol. Here we get a generalized version for $R_\alpha(a)\cap Tran$ and the studied system is more general. As a more refinement, see   Theorem  \ref{Maintheorem222-2017-levelsets} below. For   Lyapunov exponents of asymptotically additive functions,   see Section \ref{Section-Lyapunovexponents}.

    \end{Rem}

 \begin{Rem}\label{remark-example-multiplytwofunction} Let $\phi,\psi$ be two continuous functions on $X$. Define $\alpha(\mu)= {\int \phi d\mu}{\int \psi d\mu}$ which may not satisfy anyone of conditions $[A.1],[A.2],[A.3].$ Define
     $$R_{\phi,\psi}(a)=\{x\in X:\,\lim_{n\rightarrow +\infty}\frac1{n^2}\sum_{i=0}^{n-1}\phi(T^ix)\sum_{i=0}^{n-1}\psi(T^ix)=a\}.$$ Then  $R_{\phi,\psi}(a)=R_\alpha(a)$ and  thus  from item (2) of Theorem \ref{Thm--RoughResult-multifrcalanalysis}, we have
$$ h_{top} (T,R_{\phi,\psi}(a)\cap Tran  )=h_{top} (T,R_{\phi,\psi}(a)   )=\sup \{h_{\rho}  |\,\,\rho\in M (T,X)\,\,and\,\,\alpha (\rho)=a\}.$$  However, it is not sure  for
 $$I_{\phi,\psi} =\{x\in X:\,\lim_{n\rightarrow +\infty}\frac1{n^2}\sum_{i=0}^{n-1}\phi(T^ix)\sum_{i=0}^{n-1}\psi(T^ix)\text{ does not converge } \},$$ since it is unknown whether  $[A.3]$ holds so that item (1)  of Theorem \ref{Thm--RoughResult-multifrcalanalysis} can not apply. But it has positive entropy if the system has positive entropy, see Lemma \ref{Lem-IC-notempty} below.

    \end{Rem}

\subsubsection{Combination between Banach Recurrence and Multifractal Analysis }
Now we start to combine Banach recurrence with multifractal analysis to get some refined results of Theorem \ref{Thm--RoughResult-multifrcalanalysis}.



\begin{maintheorem}\label{Maintheorem111-2017-irregular}

 Suppose that $(X,\,T)$   has transitively-convex-saturated property and  entropy-dense property. Let $\alpha:M(T,X)\rightarrow \mathbb{R}$ be a continuous function  satisfying   $[A.3]$ and $Int(L_\alpha)\neq \emptyset.$
  If there is an invariant measure
with full support,
   then   $\{QR, BR_1, BR_2,BR_3, BR_4,BR_5 \}$ has   full entropy gap w.r.t. $I_\alpha\cap Tran $.

 \end{maintheorem}


\begin{Rem}\label{remark-example-multiplytwofunction-2222} For a general continuous $\alpha$ without assuming condition $[A.3],$  by slight modification in the proof  one can get that  $\{QR, BR_1, BR_2,BR_3, BR_4,BR_5 \}$ has   positive entropy gap w.r.t. $I_\alpha\cap Tran $, that is $h_{top}((BR_{i+1}\setminus BR_i)\cap I_\alpha)>0$. For example,
if $\alpha$ is   the function  as in Remark \ref{remark-example-multiplytwofunction} and assume that
 $$\inf_{\mu\in M(T,X)}{\int \phi d\mu}{\int \psi d\mu}<\sup_{\mu\in M(T,X)}{\int \phi d\mu}{\int \psi d\mu},$$ then
  in this case $$I_\alpha=I_{\phi,\psi} =\{x\in X:\,\lim_{n\rightarrow +\infty}\frac1{n^2}\sum_{i=0}^{n-1}\phi(T^ix)\sum_{i=0}^{n-1}\psi(T^ix)\text{ does not converge } \}.$$
   But it is unknown for $\{QR, BR_1, BR_2,BR_3, BR_4,BR_5 \}$ w.r.t. $R_\alpha, R_\alpha(a).$
    \end{Rem}

\begin{maintheorem}\label{Maintheorem222-2017-levelsets}
Suppose that $(X,\,T)$   has transitively-convex-saturated property and  entropy-dense property.
Let $\alpha:M(T,X)\rightarrow \mathbb{R}$ be a continuous function  satisfying  $[A.1]$, $[A.2]$ and $Int(L_\alpha)\neq \emptyset.$ Let $a\in Int(L_\alpha)$.
 If there is an invariant measure
with full support,    then
     $\{\emptyset, QR\cap BR_1 ,BR_1, BR_2,BR_3, BR_4,BR_5 \}$ has   full entropy gap w.r.t. $R_\alpha(a)\cap Tran $.


 \end{maintheorem}

\begin{maintheorem}\label{Maintheorem333-2017-alpha-regularsets}
Suppose that $(X,\,T)$   has transitively-convex-saturated property and  entropy-dense property.
Let $\alpha:M(T,X)\rightarrow \mathbb{R}$ be a continuous function  satisfying  $[A.1]$, $[A.2]$.
 If $(X,\,T)$ is not uniquely ergodic and there is an invariant measure
with full support,    then
     $\{\emptyset, QR\cap BR_1 ,BR_1, BR_2,BR_3, BR_4,BR_5 \}$ has   full entropy gap w.r.t. $R_\alpha\cap Tran $.


 \end{maintheorem}

\subsection{Applications to mixed multifractal analysis of Lyapunov exponents}\label{Section-Lyapunovexponents}


Let us recall the concept of asymptotically additive introduced from \cite{FengH}, which helps us to study multifractal behavior of Lyapunov exponents.
\begin{Def}\label{Def-asymptotically-additive}
 A sequence   of functions $\phi_n : X\rightarrow R$ is said to be asymptotically additive if
for each $\epsilon >0, $ there exists a continuous function $\phi: X\rightarrow R$ such that
$$\limsup_{n\rightarrow \infty} \frac1n \sup_{x\in X}
|\phi_n(x)-S_n\phi(x)|\leq \epsilon, $$
where $S_n\phi =  \sum_{k=0}^{n-1}\phi\circ T^k.$
\end{Def}

\subsubsection{Asymptotically Additive Functions and Lyapunov Exponents}\label{section-additiveasymptotically}
   Given a asymptotically additive sequence of functions $\Phi:=\{\phi_n:X\rightarrow \mathbb{R}\}_{n\geq 1}$,   the limit $$\lim_{n\rightarrow \infty}\frac1n\phi_n(x)$$ (if exists) is called     the Lyapunov exponent of $\Phi$ at $x$ , see \cite{FengH}.
  Define $$\alpha(\mu):=\liminf_{n\rightarrow \infty}\int\frac1n\phi_n d\mu,\,\mu\in M(T,X).$$ It is not difficult to see that for any invariant $\mu$,
  the limit $\alpha(\mu):=\lim_{n\rightarrow} \frac1n\int\phi_n(x)d\mu$ exists and the function $\alpha(\cdot):M(T,X)\rightarrow\mathbb{R}$ is affine and continuous, for example, see \cite{FengH}.
  Thus $\alpha$ satisfies conditions $[A.1]$ $[A.2]$ and $[A.3]$.  The concept of asymptotically sub-additive potentials of \cite{FengH} is mainly motivated by some
works on the Lyapunov exponents of matrix products \cite{Feng02,Feng09,FengLau} and the Lyapunov
exponents of differential maps on nonconformal repellers \cite{BarGelfert} so that the results of present paper are suitable for the cases of \cite{Feng02,Feng09,FengLau,BarGelfert}.
  In this situation,   $$R_\alpha :=\{x\in X\,| \text{ averages } \frac1n \phi_n(x)   \text{ converge as } n\rightarrow +\infty\}$$ (For convenience,  called      Lyapunov-$\Phi$-regular  set   or simply,  $\Phi$-regular   set).
   the corresponding $\Phi$-irregular set (or Lyapunov-$\Phi$-irregular set) is $I_\alpha: =X\setminus R_\alpha .$
 For any $a\in\mathbb{R},$    the level set of Lyapunov exponents is
 $$R_{\alpha}(a)  :=\{x\in X|\,\,\lim_{n\rightarrow \infty}\frac1n\phi_n(x)=a\}.$$
  There are lots of  classical results for  $R_{\alpha}(a) $, for example, see \cite{BCW,FengH}(see   \cite{PS,TV,Ol} for additive functions and a survey article \cite{Climen} and references therein) and thus our results of present paper can be as refined generalizations combined with different recurrences.



 \subsubsection{Mixed Multifractal Analysis}

 Given two asymptotically additive sequences of functions $\Phi:=\{\phi_n:X\rightarrow \mathbb{R}\}_{n\geq 1}$ and $\Psi:=\{\psi_n:X\rightarrow \mathbb{R}\}_{n\geq 1}$ where $\psi_n$ are all positive functions,   define $$\alpha(\mu):=\frac{\liminf_{n\rightarrow \infty}\int\frac1n\phi_n d\mu}{\liminf_{n\rightarrow \infty}\int\frac1n\psi_n d\mu},\,\mu\in M(T,X).$$ Similar as  discussion in subsection \ref{section-additiveasymptotically},
 $\alpha$ is a continuous function and it is not difficult to check that conditions $[A.1]$ and $[A.2]$ hold (and $[A.3]$ holds, since $[A.1]$ implies $[A.3]$). We emphasize that here $\alpha$ may be not affine.   In this case $$I_\alpha=I^\Phi_\Psi :=\{x\in X:\,\lim_{n\rightarrow +\infty} \frac{\phi_n(x)}{\psi_n(x)}\text{ does not converge } \}, R_\alpha=R^\Phi_\Psi:=X\setminus I^\Phi_\Psi\,\text{ and }$$  $$R_\alpha(a) =R^\Phi_\Psi(a):=\{x\in X:\,\lim_{n\rightarrow +\infty} \frac{\phi_n(x)}{\psi_n(x)}=a\}.$$ By Theorem \ref{Thm--PropMain-productproperty-uniformseparation},  Theorem \ref{Maintheorem111-2017-irregular},   Theorem  \ref{Maintheorem222-2017-levelsets} and  Theorem \ref{Maintheorem333-2017-alpha-regularsets},
   we have following.

\begin{maintheorem}\label{Maintheorem222-2017-Lyapunovexponents}
Suppose that $(X,\,T)$ satisfies $g$-almost product property.
  Let $\Phi:=\{\phi_n:X\rightarrow \mathbb{R}\}_{n\geq 1}$ and $\Psi:=\{\psi_n:X\rightarrow \mathbb{R}\}_{n\geq 1}$ be  two asymptotically additive sequences of functions
such that $I^\Phi_\Psi\neq \emptyset$ where $\psi_n$ are all positive functions.
   If there is an invariant measure
with full support,    then \\
    (1) $\{QR , BR_1, BR_2,BR_3, BR_4,BR_5 \}$ has   full entropy gap w.r.t. $I^\Phi_\Psi \cap Tran $. \\ 
  (2) $\{\emptyset, QR\cap BR_1 ,BR_1, BR_2,BR_3, BR_4,BR_5 \}$ has   full entropy gap w.r.t. $R^\Phi_\Psi(a)\cap Tran $ where $a\in Int(L_\alpha)$. \\
  (3) $\{\emptyset, QR\cap BR_1 ,BR_1, BR_2,BR_3, BR_4,BR_5 \}$ has   full entropy gap w.r.t. $R^\Phi_\Psi \cap Tran $.

 \end{maintheorem}

\begin{Rem}\label{Rem-Barreira-Saussol-mixed2222} If $\alpha$ is   the function  as in Remark \ref{Rem-Barreira-Saussol-mixed} for which  the part for $R_\alpha(a)$ has been discussed  in \cite{Barreira-S2001} for hyperbolic systems by Barreira and Saussol. Here we   generalize the result of \cite{Barreira-S2001} to the case of mixed   Lyapunov exponents of asymptotically additive functions and moreover various recurrences are considered simultaneously.

    \end{Rem}



\subsection{Examples}


\subsubsection{Dynamics with specification}
From  \cite{Data}, we know that for any dynamical system with specification property (not necessarily Bowen's strong version), the almost periodic points are dense in $X$  and the invariant measures supported on minimal sets are dense in the space of invariant measures.  By Lemma \ref{Lem-general--notfullsupport} (see below) there is some invariant measure with full support.
From \cite[Proposition 2.1]{PS}  we know specification implies $g$-almost product property so that the assumptions of $g$-almost product property and existence of a measure with full support in Theorem \ref{Thm--PropMain-productproperty-uniformseparation} can be replaced by specification. From \cite[Proposition 2.3 (1)]{PS2005}, entropy-dense property holds
 for dynamical systems with $g$-product property.   Thus,

\begin{maincorollary}\label{Cor-case-specification1}
 Suppose that $(X,T)$ satisfies specification  property and is not uniquely ergodic.
 Then Theorem \ref{Thm--PropMain-productproperty-uniformseparation} (1),
 Theorem \ref{Maintheorem000-2017-banachrecurrence}, Theorem \ref{Maintheorem111-2017-irregular},   Theorem  \ref{Maintheorem222-2017-levelsets} and Theorem \ref{Maintheorem333-2017-alpha-regularsets}   hold.

 \end{maincorollary}

\begin{maincorollary}\label{Cor-case-specification-plus-uniformseparation}
For non-uniquely ergodic dynamical systems with specification and uniform separation, all results of Theorem \ref{Thm--PropMain-productproperty-uniformseparation}, Theorem \ref{Maintheorem000-2017-banachrecurrence}, Theorem \ref{Maintheorem111-2017-irregular},   Theorem  \ref{Maintheorem222-2017-levelsets} and Theorem \ref{Maintheorem333-2017-alpha-regularsets}   hold.

 \end{maincorollary}

%



Corollary \ref{Cor-case-specification1} applies in following examples:
(1) It is known  from  \cite{Bl,Buzzi} that any topologically mixing interval map satisfies Bowen's  specification but  maybe not have uniform separation.
 For example, Jakobson   \cite{Jakobson} showed that there exists a set of parameter values $\Lambda\subseteq [0,4]$ of positive Lebesgue measure such
 that if $\lambda\in \Lambda,$ then the logistic map $f_\lambda(x)=\lambda x(1-x)$ is topologically mixing.
 (2) Recently we learned from \cite{GL} that $C^0$ generic volume-preserving dynamical systems have specification property.
  (3)
In \cite{Aoki}  N. Aoki constructs a zero-dimensional
ergodic automorphism without densely periodic property that  obeys   specification, but not Bowen's specification.
For the class of all solenoidal automorphisms, it is proved in \cite{AokiDK}
that the class of automorphisms with  specification is wider than
the class of automorphisms with Bowen's specification. 
(4) In \cite{PS}, Pfister and   Sullivan  gave an example of
a dynamical system with finite topological entropy, for which the entropy density of ergodic
measures is true (specification property is true), but the uniform separation property and
the upper semi-continuity of the entropy map fail. This example is a subshift of the
shift space $Y := [-1, 1]^{\mathbb{Z}^+}$, see \cite[Page 952-Page 953]{PS} for more details.

Corollary \ref{Cor-case-specification-plus-uniformseparation} applies in following examples:  (1)
  all topological mixing subshifts of finite type and mixing soft subshifts (in particular,  all full shifts on finite alphabets);
  (2)
   all  subsystems restricted on topological mixing locally maximal expanding set or hyperbolic set (in particular, all topological mixing  expanding maps or topological mixing hyperbolic diffeomorphisms, called Anosov).
  This is because such examples all satisfy  Bowen's  specification (see \cite{Bowen70,Weiss} for item (1) and   see \cite{Bowen71-trans} for item (2), also cf. \cite{DGS}) and are expansive  which is stronger than uniform separation, see \cite{PS}.
(3)
 The results of this paper are also applicable to some dynamical systems beyond uniform hyperbolicity.
From \cite{Go} we know that  non-hyperbolic diffemorphism $f$ with $C^{1+Lip}$ smoothness, conjugated to
 a transitive Anosov diffeomorphism $g$, exists even the conjugation and its inverse  is H$\ddot{\text{o}}$lder continuous.
(4)
  From   \cite[Section 4.3]{To2010} we know the time-1 map of a transitive Anosov flow
     satisfies specification. 
     In this case, $f$ is partially hyperbolic with one-dimension central bundle. Then $f$ is far from tangency so that  $f$ is entropy-expansive   (see  \cite{LiaoVianaYang} or see    \cite{DFPV,PacVie}). Recall that from   \cite{Misiurewicz} entropy-expansive implies  asymptotically $h-$expansive and from
      Theorem 3.1  of    \cite{PS} any expansive or asymptotically
      $h-$expansive system satisfies uniform separation property. 




\subsubsection{$\beta$-shifts}

    \begin{maincorollary}\label{Cor-main-beta}
  Theorem \ref{Thm--PropMain-productproperty-uniformseparation}, Theorem \ref{Maintheorem000-2017-banachrecurrence}, Theorem \ref{Maintheorem111-2017-irregular},   Theorem  \ref{Maintheorem222-2017-levelsets} and Theorem \ref{Maintheorem333-2017-alpha-regularsets}   hold for any $\beta-$shfit.
 \end{maincorollary}

 Let us explain why  Corollary \ref{Cor-main-beta} holds. We know that any $\beta-$shfit
    is expansive (which is stronger than uniform separation, see \cite{PS}) and from  \cite{PS}   it always satisfies $g$-almost product property.
   Furthermore, from \cite{Sig-76}  we know that periodic points are dense in the whole space and the periodic measures are dense in the space of invariant measures (i.e., $\overline{ M_{p} (T,X) }=M (T,X)).$ By Lemma \ref{Lem-general--notfullsupport} (see below) there is some invariant measure with full support (from \cite [Theorem 13 (ii)]{Walters78} we can also learn that  the unique maximal entropy measure  of $\beta-$shifts always carries full support). Thus  Theorem \ref{Thm--PropMain-productproperty-uniformseparation}, Theorem \ref{Maintheorem000-2017-banachrecurrence}, Theorem \ref{Maintheorem111-2017-irregular},   Theorem  \ref{Maintheorem222-2017-levelsets} and Theorem \ref{Maintheorem333-2017-alpha-regularsets}   apply  in any $\beta-$shfit.
  It is worth mentioning that from    \cite{Buzzi} the set of parameters of $\beta$ for which   specification holds, is dense in $ (1,+\infty)$ but has Lebesgue zero measure.

  It is well-known that any subshift of finite type and any $\beta-$shfit have unique maximal entropy measure with full support. From \cite{Pavlov} there exist  subshifts  with
multiple measures of maximal entropy with disjoint support
which have $g$-almost product property and a measure with full support. Thus    Theorem \ref{Thm--PropMain-productproperty-uniformseparation}, Theorem \ref{Maintheorem000-2017-banachrecurrence}, Theorem \ref{Maintheorem111-2017-irregular},   Theorem  \ref{Maintheorem222-2017-levelsets} and Theorem \ref{Maintheorem333-2017-alpha-regularsets} also apply  in these subshifts \cite{Pavlov}.

\subsection{Organization of this paper} In Section \ref{section-prepare}  we recall some notions, lemmas and in  Section \ref{section-transitivelysaturated} we give the proof of  Theorem \ref{Thm--PropMain-productproperty-uniformseparation}. In Section \ref{section-prove-banach-multifractal} we will study Banach recurrent gap-sets w.r.t. irregular set and level set respectively. In Section \ref{section-commentsquestions} we give some comments and further questions.

\section{Preliminaries}\label{section-prepare}

\subsection{Invariant measures}
Let ${\{\varphi_{j}\}}_{j\in\mathbb{N}}$ be a dense subset of $C(X,\mathbb{R})$ which is the space of continuous functions, then
$$d (\xi,\tau)=\sum_{j=1}^{\infty}\frac{|\int\varphi_{j}d\xi-\int\varphi_{j}d\tau|}{2^{j+1}\|\varphi_{j}\|}$$
defines a metric on $M(X)$ for the $weak^{*}$ topology  \cite{Walter}, where
 $\|\varphi_i\|=\max\{|\varphi_i(x)|:x\in X\}.$
Note that
\begin{equation}\label{diameter-of-Borel-pro-meas}
d (\xi,\tau)\leq 1~~\textrm{for any}~~\xi,\tau\in M(X).
\end{equation}


\subsection{Entropy}


Let $T:X\rightarrow X$  be a continuous map of a compact
metric space $X$.
Now let us to recall the definition of topological entropy  in  \cite{Bowen1} by Bowen.
Let $x\in X$. The dynamical ball $B_n (x,\varepsilon)$ is the set $$B_n (x,\varepsilon):=\{y\in X|\,\max\{d (T^j (x),T^j (y))|\,0\leq j\leq n-1\}\leq \varepsilon\}.$$
Let $E\subseteq X,$ and $\mathfrak{F}_n (E,\epsilon)$ be the collection of all finite or countable covers of $E$ by sets of the form $B_m (x,\epsilon)$ with $m\geq n.$ We set $$C (E;t,n,\epsilon,T):=\inf\{\sum_{B_m (x,\epsilon)\in \mathcal{C}}2^{-tm}:\mathcal{C}\in \mathfrak{F}_n (E,\epsilon)\},\,\,\,\,\,\text{ and }$$    $ C (E;t,\epsilon,T):=\lim_{n\rightarrow \infty}C (E;t,n,\epsilon,T).$  Then  $h_{top} (E,\epsilon,T):=\inf\{t:C (E;t,\epsilon,T)=0\}=\sup\{t: C (E;t,\epsilon,T)=\infty\}$  and the {\it topological entropy} of $E$ is defined as $$h_{top} (T,E);=\lim_{\epsilon\rightarrow 0}h_{top} (E,\epsilon,T).$$ In particular, if $E=X,$ we also denote $h_{top} (T,X)$ by $h_{top} (T)$.
It is known   from    \cite{Bowen1} that if $E$ is an invariant compact subset, then the topological entropy $h_{top} (T,E)$ is same as the classical definition (for classical definition of topological entropy, see Chapter 7 in \cite{Walter}). Let us recall some basic facts about topological entropy.
From    \cite{Bowen1} we know for any $Y\subseteq X,$ \begin{eqnarray}\label{Bowen-fY=Y}
h_{top} (T,f Y)= h_{top} (T,Y).
\end{eqnarray}and
for any subsets $Y_1\subseteq Y_2\subseteq X,$
\begin{eqnarray}\label{Y1containY2}
h_{top} (T,Y_1)\leq h_{top} (T,Y_2).
\end{eqnarray}
 If one considers a collection of subsets of $X$:   $\{Y_i\}_{i=1}^ {+\infty},$    from   \cite{Bowen1} we know that    the topological entropy satisfies \begin{eqnarray}\label{Bowen-disjointSet-entropy}
 h_{top} (T,\bigcup_{i=1}^{+\infty}Y_i)=\sup_{ i\geq 1} h_{top} (T,Y_i).
 \end{eqnarray}

 Let $\xi=\{V_i|\,i=1,2,\cdots,k\},$ be a finite partition of measurable sets of $X$. The entropy of a probability measure $\nu\in M(X)$ with respect to $\xi$
 is  $H(\nu,\xi):=-\sum_{V_i\in\xi}\nu(V_i)\log \nu(V_i).$
 We write $T^{\vee n}\xi:=\vee_{k\in \Lambda }T^{-k}\xi.$ The entropy of an invariant measure $\nu\in M(T,X)$ with respect to $\xi$ is $h(T,\nu,\xi):=\lim_{n\rightarrow \infty}\frac 1n H(\nu, T^{\vee n}\xi),$  and the {\it metric entropy} of $\nu$ is $$h_\nu(T):=\sup_{\xi} h(T,\nu,\xi).$$
 More information of metric entropy, see Chapter 4 of  \cite{Walter}.

 For convenience,  we write $h_{top}$ and $h_\omega$ to denote $h_{top}(T)$ and $h_\omega(T)$.



\subsection{Specification property and Product property}

Firstly we recall  the definition of   specification property which is stronger than  $g-$almost product property, see \cite{DGS,Sig,Bow,Bowen2,Bowen71-trans,To2010}. Let  $\,T$ be a      continuous map of a compact metric space $X$.

\begin{Def}\label{specification}  We say that the dynamical system $T$  satisfies {\it specification property}, if the following holds:  for any $\epsilon>0$ there exists an integer $M(\epsilon)$ such that for any $k
\geq 2,$ any $k$ points $x_1,\cdots,x_k$, any  integers $$a_1\leq b_1<a_2\leq b_2\cdots<a_k\leq b_k$$ with $a_{i+1}-b_i\geq M(\epsilon)\,(1\leq i\leq k-1),$   there exists a point $x\in X$ such that \begin{eqnarray}\label{specification-inequality}
 d(T^j(x),T^j(x_i))
 \leq \epsilon,\,\,\,\, \text{ for } \,\,a_i\leq j\leq b_i,\,1\leq i\leq k.
 \end{eqnarray}
\begin{eqnarray}\label{specification-inequality-Bowenball-language}
\text{ In other words, the set }  \hat{B}=\cap_{i=1}^k f^{-a_i} B_{b_i-a_i}(f^{a_i} x_i,\epsilon)  \text{ is nonempty. }
 \end{eqnarray}

\end{Def}

The original definition of specification, due to Bowen, was stronger.

\begin{Def}\label{Bowen-specification} We say that the dynamical system $T$  satisfies {\it Bowen's  specification property}, if  under the assumptions of Definition \ref{specification} and   for any  integer $p\geq M(\epsilon)+b_k-a_1,$ there exists a point $x\in X$ with $T^p(x)=x$ satisfying (\ref{specification-inequality}).

\end{Def}

Now we start to recall the concept $g-$almost product property in  \cite{PS}  (there is a slightly weaker variant,  called almost specification, see  \cite{Tho2012}).  It is weaker than specification property (see  Proposition 2.1 in \cite{PS}).  A striking and typical example of   $g-$almost product property (and almost specification) is that it applies to every $\beta-$shift  \cite{PS,Tho2012}. In sharp contrast, the set of $\beta$ for which
the $\beta-$shift has specification property has zero Lebesgue measure   \cite{Buzzi,Schmeling}.

Let $\Lambda_n=\{0,1,2,\cdots,n-1\}.$ The cardinality of a finite set $\Lambda$ is denoted by $\# \Lambda.$ Let $x\in X$. The dynamical ball $B_n(x,\varepsilon)$ is the set $$B_n(x,\varepsilon):=\{y\in X|\,\max\{d(T^j(x),T^j(y))|\,j\in\Lambda_n\}\leq \varepsilon\}.$$

\begin{Def}\label{blowup-function} Let $g:\mathbb{N}\rightarrow \mathbb{N}$  be a given nondecreasing unbounded map with the properties $$g(n)<n\,\,\text{ and } \lim_{n\rightarrow \infty}\frac{g(n)}n=0.$$ The function $g$ is called {\it blowup function.} Let $x\in X$ and $\varepsilon>0.$ The $g-$blowup of $B_n(x,\varepsilon)$ is the closed set
 $B_n(g;x,\varepsilon):=$ $$\{y\in X|\, \exists \Lambda\subseteq\Lambda_n ,\#(\Lambda_n\setminus\Lambda)\leq g(n)\,\text{ and }\,\max\{d(T^j(x),T^j(y))|\,j\in\Lambda\}\leq \varepsilon\}.$$

\end{Def}

\begin{Def}\label{product-property} We say that the dynamical system $T$  satisfies {\it $g-$almost product property}  with blowup function $g$, if  there is a nonincreasing function $m:\mathbb{R}^+\rightarrow \mathbb{N},$ such that
for any $k
\geq 2,$ any $k$ points $x_1,\cdots,x_k\in X$, any positive $\varepsilon_1,\cdots,\varepsilon_k$ and any integers $n_1\geq m(\varepsilon_1),\cdots, n_k\geq m(\varepsilon_k),$ $$\bigcap_{j=1}^k T^{-M_{j-1}}B_{n_j}(g;x_j,\varepsilon_j)\neq \emptyset,$$ where $M_0:=0,M_i:=n_1+\cdots+n_i,i=1,2,\cdots,k-1.$

\end{Def}

It is  well known that the natural projection $x\mapsto \delta_x$ is continuous and if we define operator $T_f$
on $M(X)$ by formula $T_f(\mu)(A)=\mu(f^{-1}(A))$, then we can identify $(X,f)$
with $T_f$ restricted to the set of Dirac measures (these systems are conjugate). Therefore, without loss of generality
we will assume that $d(x,y)=d (\delta_x,\delta_y)$. Denote  a ball in $M(X)$ by
$$\mathcal B(\nu,\zeta):=\{\mu\in M(X):d(\nu,\mu)\leq \zeta\}.$$ For $x\in X,$ define  $\Upsilon_n (x):=\frac1n\sum_{j=0}^{n-1}\delta_{T^j (x)}$  where
  $\delta_y$ is the Dirac  probability measure supported
  at $y\in X$.   

\begin{lem}\label{lem:Approx-by-OrbitMeasure}
	
\cite[Lemma 2.1]{PS}
Assume that $(X,T)$ satisfies $g$-almost product property. Let
$x_1,\cdots,x_k\in X$,    $\varepsilon_1,\cdots,\varepsilon_k$ and
 $q_1\geq m(\varepsilon_1),\cdots, q_k\geq m(\varepsilon_k)$
 be given. Assume that
 $$\Upsilon_{q_j}(x_j)\in \mathcal B(\nu_j,\zeta_j),\,j=1,2,\cdots,k. $$
 Then for any $z\in \cap_{j=1}^k T^{-Q_{j-1}} B_{q_j}(g;x_j,\epsilon_j)$ and any probability measure $\alpha$
 $$d(\Upsilon_{Q_k}(z),\alpha)\leq \sum_{j=1}^k \frac{n_j}{Q_k}(\zeta_i+\epsilon_i+\frac{g(q_i)}{q_i})$$
 where $Q_0=0,Q_i=q_1+\cdots+q_i.$

\end{lem}

\subsection{Uniform  separation}

Now we recall the definition of uniform separation property  \cite{PS}.
For $\delta>0$ and  $\varepsilon>0$, two points $x$ and $y$ are
$(\delta,n,\varepsilon)-$separated if $$\#\{j:d(T^jx,T^jy)>\varepsilon,\,j\in\Lambda_n\} \geq \delta n.$$ A subset
$E$ is  $(\delta,n,\varepsilon)$-separated if any pair of different points of $E$ are  $(\delta,n,\varepsilon)-$separated.      Let $F\subseteq M(X)$ be a neighborhood of $\nu\in M(T,X)$.
 Define $$ X_{n,F}:=\{x\in X|\, \Upsilon_n(x)\in F\},$$
and define
 \begin{eqnarray*}N(F;\delta,n,\varepsilon):=\text{maximal cardinality of a } (\delta, n,\varepsilon)-\text{separated subset of } X_{n,F}.\end{eqnarray*}

\begin{Def}\label{uniform-separation-property} We say that the dynamical system $T$  satisfies {\it uniform separation  property}, if  following holds. For any $\eta>0,$ there exist $\delta^*>0,\epsilon^*>0$ such  that for $\mu$ ergodic and any neighborhood $F\subseteq M(X)$ of $\mu$, there exists $n^*_{F,\mu,\eta},$ such that for $n\geq n^*_{F,\mu,\eta},$ $$N(F;\delta^*,n,\epsilon^*)\geq 2^{n(h_\mu(f)-\eta)}.$$

\end{Def}

\begin{Lem}\label{Lemm-PS-uniformseparation}\cite[Corollary 3.1]{PS} Assume that $(X,d)$ has the uniform separation property, and has entropy-dense property. For any $\eta,$ there exist $\delta^*>0,\epsilon^*>0$ such  that for $\mu\in M(T,X)$   and any neighborhood $F\subseteq M(X)$ of $\mu$, there exists $n^*_{F,\mu,\eta},$ such that for $n\geq n^*_{F,\mu,\eta},$ $$N(F;\delta^*,n,\epsilon^*)\geq 2^{n(h_\mu(f)-\eta)}.$$

\end{Lem}

\section{Proof of Theorem \ref{Thm--PropMain-productproperty-uniformseparation}: Minimal Set and Saturated Property}\label{section-transitivelysaturated}

\subsection{Full support and minimal set}

\begin{Lem}\label{Lem-saturated-imply-transitive}    Suppose that $T$ is single-saturated and there is some invariant measure $\mu$ with full support (i.e., $S_\mu=X $).
 Then  $T$ is transitive.

\end{Lem}

{\bf Proof.} By single-saturated property, $G^T_\mu\neq \emptyset$. Take $x\in G^T_\mu$, then $\omega_T(x)=X.$\qed  

 \begin{Lem}\label{Lem-general--notfullsupport}    Suppose that a subset $B'$ of  $B:=\{\omega\in M(T,X)|\,S_\omega \neq X \}$
 is dense in $M(T,X)$.  Then there is some invariant measure $\mu$ with full support (i.e., $S_\mu=X $) $\Leftrightarrow$ $\overline{\cup_{\omega\in B'}S_\omega}=X.$

\end{Lem}
{\bf Proof.}  $\Rightarrow$ By assumption there is a sequence of invariant measures $\mu_i\in B'$ with  $S_{\mu_i}\neq X$ converging to $\mu$. Then
$1=\limsup_{n\rightarrow}\mu_n(\overline{\cup_{\omega\in B'}S_\omega})\leq \mu(\overline{\cup_{\omega\in B'}S_\omega}).$ It follows that $X=S_\mu\subseteq \overline{\cup_{\omega\in B'}S_\omega}$.

$\Leftarrow$ Take a sequence of invariant measures $\mu_i\in B'$
 with  $S_{\mu_i}\neq X$ such that $ \overline{\cup_{i\geq 1}S_{\mu_i}}=X.$ Let $\mu=\sum_{n\geq 1}\frac1{2^n}\mu_n.$ Then $\mu(\cup_{i\geq 1}S_{\mu_i})=1$ so that $S_\mu=X.$ \qed

 \begin{Rem} Similarly, one can get that for any dynamical system $(X,T),$ there is an invariant measure with full support $\Leftrightarrow$ $\overline{\cup_{\mu\in M(T,X)}S_\mu}=X.$

 \end{Rem}


\begin{Lem}\label{Lem-notfullsupport}    Suppose that $T$ has entropy-dense property and there is some invariant measure $\mu$ with full support (i.e., $S_\mu=X $).
Then $B:=\{\omega\in M_{erg}(T,X)|\,S_\omega \neq X \}$
  is dense in $M(T,X)$ and $\overline{\cup_{\omega\in B}S_\omega}=X.$

\end{Lem}

{\bf Proof.}  By entropy-dense property, the set $B$ is dense in $M(T,X).$ Then by Lemma \ref{Lem-general--notfullsupport} $\overline{\cup_{\omega\in B}S_\omega}=X.$ \qed

\begin{Lem}\label{Lem-minimalpointsdense}    Suppose that    $\{\mu\in M(T,X)|\,S_\mu \text{ is minimal }\}$ is dense in $M(T,X)$. Then there is some invariant measure $\mu$ with full support (i.e., $S_\mu=X $) $\Leftrightarrow$  almost periodic set $A$ is dense in $X$.

\end{Lem}

{\bf Proof.}  $\Rightarrow$ By assumption there is a sequence of invariant measures $\mu_i$  with $S_{\mu_i}\subseteq A$    converging to $\mu$. Then
$1=\limsup_{n\rightarrow}\mu_n(A)\leq \limsup_{n\rightarrow}\mu_n(\overline{A})\leq\mu(\overline{A}).$ It follows that $X=S_\mu\subseteq \overline{A}$.

$\Leftarrow$ Take a sequence of points $\{x_i\}\subseteq A$ dense in $X$. For any $i$, take $\mu_i$ to be a invariant measure on $\omega_T(x_i).$ Then    $x_i\in\omega_T(x_i)=S_{\mu_i} $ and so $ \overline{\cup_{i\geq 1}S_{\mu_i}}=X.$ Let $\mu=\sum_{n\geq 1}\frac1{2^n}\mu_n.$ Then $\mu(\cup_{i\geq 1}S_{\mu_i})=1$ so that $S_\mu=X.$\qed

\begin{Lem}\label{Lem-minimalmeasuredense}     Suppose that $T$ has $g$-almost product property.
Then ergodic measures supported on minimal sets are dense in $M(T,X)$.

\end{Lem}

{\bf Proof.}  Let  $\nu\in M (T,X)$ and  $G\subseteq M (X)$ be a neighborhood  of $\nu.$
 Take an open ball  $G'\subseteq M(X)$ such that
 $\nu\in G'\subseteq  \overline{G'}\subset G$.
 From the proof of       \cite[Proposition 2.3 (1)]{PS2005}, one construct  a closed invariant set $Y$ and there exists $n_{G'} \in\mathbb{N}$ such that   for
any $y\in Y$ and any $n\geq n_{G'},$    $\Upsilon_n (y)\in G'.$ So for
 any $m\in M_{erg} (T,Y),$ by Birkhoff ergodic theorem there is $y\in Y$ such that $ \Upsilon_n (y)$ converge to $m$ in weak$^*$ topology and thus $m\in \overline{G'}.$ In other words, $M_{erg} (T,Y)\subseteq \overline{ G'}.$  By convex property of the ball $G'$ and Ergodic Decomposition theorem,  $M (T,Y)\subseteq \overline{ G'}.$
 Take an ergodic measure $\mu$ supported on a minimal subset of $Y$, then $\mu\in \overline{ G'}\subseteq G.$
\qed

\begin{Lem}\label{Lem-minimalpoints-measuresdense}    Suppose that $T$ has $g$-almost product property and there is some invariant measure $\mu$ with full support (i.e., $S_\mu=X $).
Then ergodic measures supported on minimal sets are dense in $M(T,X)$ and almost periodic set $A$ is dense in $X$.

\end{Lem}
{\bf Proof.} by Lemma \ref{Lem-minimalmeasuredense} ergodic measures supported on minimal sets are dense in $M(T,X)$. Combining with Lemma \ref{Lem-minimalpointsdense},  almost periodic set $A$ is dense in $X$.\qed

\subsection{Transitively saturated: Proof of Theorem \ref{Thm--PropMain-productproperty-uniformseparation}}\label{Section-Proposition-MAIN}

 By Lemma \ref{Lem-minimalpoints-measuresdense} the set  $$\{\omega\in M(T,X)|\,\mu\text{ is ergodic, } \,S_\omega\text{ is minimal }\}$$
  is dense in $M(T,X)$  and almost periodic set $A$ is dense in $X$. 
  So we only need to prove following Theorem \ref{Thm-transtively-saturated} and Theorem \ref{Thm-transtively-convex-saturated} which imply Theorem \ref{Thm--PropMain-productproperty-uniformseparation}.

 \begin{Thm}\label{Thm-transtively-saturated}
  Let $T:X\rightarrow X$ be  a continuous map of a compact metric space  $X$ with  $g$-almost product property  property and uniform separation.
 Suppose that there is an invariant measure with full support.     Then $T$ is transitively-saturated.

\end{Thm}

{\bf Proof}. From  \cite[Theorem 4.1 (3)]{PS} $h_{top} (f,G_K)\leq \inf\{h_\mu (f)\,|\,\mu\in K\}$.  Since  $G^T_K $ is contained
in $G_K,$  then $h_{top} (f,G^T_K)\leq h_{top} (f,G_K)\leq\inf\{h_\mu (f)\,|\,\mu\in K\}$.

The difficult part of the proof is to obtain a lower bound for $h_{top}(T ,G^{T}_K).$
One can modify the construction in the proof of \cite[Theorem 1.1]{PS} to obtain a closed subset $F\subseteq G_K$ such that
the entropy of $F$ close to $\inf\{h_\mu (f)\,|\,\mu\in K\}$ and simultaneously we can require that the chosen points in $F$ is transitive.
For convenience of readers, we give a precise construction as follows.

By connectedness of $K$, one has

\begin{Lem}\cite[Page 944]{PS} (or \cite[Page 202]{DGS}) There exists a sequence $\{\alpha_1,\alpha_2,\cdots\}$ in $K$  such that
$$\overline{\{\alpha_j: j\in \mathbb{N}, \,j>n\}}=K, \,\,\forall \,n\in \mathbb{N}\,\text{ and } \lim_{j\rightarrow\infty} d(\alpha_j,\alpha_{j+1})=0.$$

\end{Lem}

Let $\eta>0$ and $$h^*:=\inf\{h_\mu (f)\,|\,\mu\in K\}-2\eta,\,\, H^*:= \inf\{h_\mu (f)\,|\,\mu\in K\}-\eta.$$ Given this sequence of measures $\{\alpha_k\}$, we will construct a subset $G$ such that for each $x\in G,$ $\{\Upsilon_n(x)\}$ has the same limit-point set as the sequence $\{\alpha_k\},$ and $h_{top}(T,G)\geq h^*.$ The construction of $G$ is the core of the proof which is also used in the proof of Theorem \ref{Thm-transtively-convex-saturated} below.

 By Lemma \ref{Lemm-PS-uniformseparation},  we can find   $\delta^*>0,\epsilon^*>0$ such  that for $\mu\in M(T,X)$  and any neighborhood $F\subseteq M(X)$ of $\mu$, there exists $n^*_{F,\mu,\eta},$ such that for $n\geq n^*_{F,\mu,\eta},$ \begin{eqnarray}\label{Eq-uniformseparation}
N(F;\delta^*,n,\epsilon^*)\geq 2^{n(h_\mu(f)-\eta)}.
\end{eqnarray}

 Let $m:\mathbb{R}^+\rightarrow \mathbb{N}$ be the nonincreasing function by $g-$almost product property.
Let $\{\zeta_k\}$ and $\{\epsilon_k\}$ be two strictly decreasing sequences so that $\lim_{k\rightarrow\infty}\zeta_k=0=\lim_{k\rightarrow\infty}\epsilon_k$ with $\epsilon_1<\frac14\epsilon^*.$ By Lemma \ref{Lem-minimalpoints-measuresdense}    almost periodic set $A$ is dense in $X$. Thus for any fixed $k$ there is a finite set $\Delta_k:=\{x^k_1,x^k_2,\cdots,x^k_{t_k}\}\subseteq A$  and $L_k\in \mathbb{N}$ such that $\Delta_k$ is $\epsilon_k-$dense in $X$ and for any $1\leq i\leq t_k,$ any $l\geq 1,$ there is $n\in[l,l+L_k]$ such that $f^n(x^k_i)\in B(x^k_i,\epsilon_k).$  This implies that any $1\leq i\leq t_k,$
\begin{eqnarray}\label{Eq-minmialrole-Mainproof}
\frac{\#\{0\leq n \leq lL_k:\, d(T^n x^k_i, x^k_i)\leq \epsilon_k\}}{lL_k}\geq \frac1{L_k}.
\end{eqnarray}
Take $l_k$ large enough such that
\begin{eqnarray}\label{Eq-minmialrole-Mainproof2222222222222}
l_kL_k\geq m(\epsilon_k), \,\,\frac{g(l_kL_k)}{l_kL_k}< \frac1{4L_k}.
\end{eqnarray}
We may assume the two sequences of $\{t_k\},\{l_k\},\{L_k\}$ are strictly increasing.

From (\ref{Eq-uniformseparation}) we get the existence of $n_k$ and a $(\delta^*, n_k, \epsilon^*)-$separated subset $\Gamma_k\subseteq X_{n_k,{\mathcal B}(\alpha_k,\zeta_k)}$ with
\begin{eqnarray}\label{Eq-Gamma_k}
\# \Gamma_k \geq 2^{n_k H^*}
\end{eqnarray}
We may assume that $n_k$ satisfies
\begin{eqnarray}\label{Eq-choose-n_k}
n_k>m(\epsilon_k),\, \frac{t_kl_kL_k}{n_k}\leq \zeta_k,\, \delta^* n_k> 2 g(n_k)+1 \text{ and } \frac{g(n_k)}{n_k}\leq \epsilon_k \,\,\,\text{ and }
\end{eqnarray}
   \begin{eqnarray}\label{Eq-choose-n_k2222}
2^{H^*n_k }
 \geq 2^{h^*(n_k+t_kl_kL_k)}.\end{eqnarray}


We choose a strictly increasing $\{N_k\}$, with $N_k\in \mathbb{N},$ so that
\begin{eqnarray}\label{Eq-choose-N_k}
   n_{k+1}+\,t_{k+1}l_{k+1} L_{k+1} \leq \zeta_k \sum_{j=1}^{k} (n_jN_j+t_jl_jL_j) \,\,\,\text{ and }
\end{eqnarray}
 \begin{eqnarray}\label{Eq-choose-N_k2222222}
  \sum_{j=1}^{k-1} (n_jN_j+t_jl_jL_j) \leq \zeta_k \sum_{j=1}^{k} (n_jN_j+t_jl_jL_j) .
\end{eqnarray}
Now we define the sequences $\{n'_j\}$,  $\{\epsilon'_j\}$ and  $\{\Gamma'_j\}$, by setting for
$$j=N_1+N_2+\cdots +N_{k-1} +t_1 +\cdots+t_{k-1} +q \text{ with }1\leq q \leq N_k,$$
$$n'_j:=n_k,\,\, \epsilon'_j:=\epsilon_k, \,\,\Gamma'_j:=\Gamma_k\,\,\,\text{ and for }$$ 
$$j=N_1+N_2+\cdots +N_{k} +t_1 +\cdots+t_{k-1} +q \text{ with }1\leq q \leq t_k,$$
$$n'_j:=l_kL_k,\,\, \epsilon'_j:=\epsilon_k, \,\,\Gamma'_j:=\{x^k_q\}.$$
 Let $$G_k:=\bigcap_{j=1}^k(\bigcup_{x_j\in \Gamma'_j} T^{-M_{j-1}} B_{n'_j}(g; x_j,\epsilon'_j))\text{ with } M_j:=\sum_{l=1}^j n'_l.$$  Note that $G_k$ is non-empty closed set. We can label each set obtained by developing this formula by the branches
of a labeled tree of height $k$. A branch is labeled by $(x_1,\cdots,x_k)$
with $x_j\in \Gamma'_j.$ Then Theorem \ref{Thm-transtively-saturated} can be deduced by following lemma.

\begin{Lem} Let $\epsilon=\frac14 \epsilon^*$ and let $G:=\cap_{k\geq 1} G_k.$
Then we have the following. \\
(1) Let $x_j,y_j\in \Gamma'_j$ with $x_j\neq y_j.$ If $x\in B_{n'_j}(g; x_j,\epsilon'_j)$ and $y\in B_{n'_j}(g;y_j,\epsilon'_j),$ then $$\max\{d(T^mx,T^my): m=0,\cdots,n_j-1\}>2\epsilon.$$
(2) $G$ is closed set that is the disjoint union of non-empty closed sets $G(x_1,x_2,\cdots)$ Labeled by
$(x_1,x_2,\cdots)$ with $x_j \in \Gamma'_j$. Two different sequences label two different sets.\\
(3) $G\subseteq G_K.$\\
(4) $h_{top}(T,G)\geq H^*-\eta=h^*$.\\
(5) $G\subseteq Tran.$

\end{Lem}

{\bf Proof.} Different with \cite[Lemma 5.1]{PS}, our new construction implies item (5).
We can modify the proof  of \cite[Lemma 5.1]{PS} adaptable to our new construction   and simultaneously the new construction
guarantees item (5).

(1) Let $x\in B_{n'_j}(g;x_j,\epsilon'_j)$ and $y\in B_{n'_j}(g;y_j,\epsilon'_j)$. Since
$x_j$ and $y_j$ are $(\delta^*,n'_j,\epsilon^*)$-separated and (\ref{Eq-choose-n_k}) holds, there exists
$m\in \Lambda_{n'_j}$ so that
$$d(T^mx_j,T^my_j)>\epsilon^*=4\epsilon,\,d(T^mx_j,T^mx)\leq \epsilon'_j,\,d(T^my_j,T^my)\leq \epsilon'_j.$$
However,
$$d(T^mx,T^my)\geq d(T^mx_j,T^my_j)-d(T^mx_j,T^mx)-d(T^my_j,T^my)>2\epsilon.$$

(2) Note that $G$ is the intersection of closed sets. Let $(x_1,x_2,\cdots)$ be a sequence with
 $x_j\in \Gamma'_j.$ By the $g-$almost product property and compactness
 $$\bigcap_{j\geq 1}  T^{-M_{j-1}} B_{n'_j}(g; x_j,\epsilon'_j)$$ is nonempty and closed. By item (1) the two sets of
$B_{n'_j}(g; x_j,\epsilon'_j)$ and $B_{n'_j}(g; y_j,\epsilon'_j)$ are
disjoint when $x_j\neq y_j.$ So two different sequences label two different sets.



(3) Define the stretched sequence
$\{\alpha'_m\}$ by
$$\alpha'_m:=\alpha_k\,\,\text{ if }\,\sum_{j=1}^{k-1} (n_jN_j+t_jl_jL_j)+1\leq m \leq \sum_{j=1}^{k} (n_jN_j+t_jl_jL_j). $$
 Then the sequence $\{\alpha'_m\}$ has the same limit-point set as the sequence of $\{\alpha_k\}.$ If
 $$\lim_{n\rightarrow\infty} d(\Upsilon_n(y),\alpha'_n)=0,$$ then the two sequences $\{\Upsilon_n(y)\},\{\alpha'_n\}$
 have the same limit-point set. By (\ref{Eq-choose-N_k})  $\lim_{n\rightarrow \infty }\frac {M_{n+1}}{M_n}=1$. So from the definition of $\{ \alpha'_n\}$, we only need to prove that
 for any $y\in G,$ one has $$\lim_{l\rightarrow\infty} d(\Upsilon_{M_l}(y),\alpha'_{M_l})=0.$$

Assume that $\sum_{j=1}^{k} (n_jN_j+t_jl_jL_j)+1\leq M_l \leq \sum_{j=1}^{k+1} (n_jN_j+t_jl_jL_j).$
If $M_l\leq \sum_{j=1}^{k} (n_jN_j+t_jl_jL_j)+n_{k+1}N_{k+1},$ by Lemma \ref{lem:Approx-by-OrbitMeasure}  and
(\ref{Eq-choose-n_k})
$$
d(\Upsilon_{  M_l- \sum_{j=1}^{k} (n_jN_j+t_jl_jL_j)}(T^{   \sum_{j=1}^{k} (n_jN_j+t_jl_jL_j)}y),\alpha_{k+1})\leq \zeta_{k+1}+2\epsilon_{k+1}.
$$
Otherwise,  $M_l>  \sum_{j=1}^{k} (n_jN_j+t_jl_jL_j)+n_{k+1}N_{k+1},$ by Lemma \ref{lem:Approx-by-OrbitMeasure},
(\ref{Eq-choose-n_k}) and  (\ref{diameter-of-Borel-pro-meas})  
\begin{eqnarray*}
& &d(\Upsilon_{  M_l- \sum_{j=1}^{k} (n_jN_j+t_jl_jL_j)}(T^{   \sum_{j=1}^{k} (n_jN_j+t_jl_jL_j)}y),\alpha_{k+1})
\\
&
\leq& \frac { n_{k+1}N_{k+1}}{  M_l- \sum_{j=1}^{k} (n_jN_j+t_jl_jL_j)}d(\Upsilon_{  n_{k+1}N_{k+1}}(T^{   \sum_{j=1}^{k} (n_jN_j+t_jl_jL_j)}y),\alpha_{k+1})
\\
& &
+\frac {M_l- \sum_{j=1}^{k} (n_jN_j+t_jl_jL_j)- n_{k+1}N_{k+1}}{  M_l- \sum_{j=1}^{k} (n_jN_j+t_jl_jL_j)}\times 1
\\
&
\leq&   1\times(\zeta_{k+1}+2\epsilon_{k+1})+\frac{t_{k+1}l_{k+1} L_{k+1}}{  n_{k+1}N_{k+1}} \\
&\leq & 2\zeta_{k+1}+2\epsilon_{k+1}.
\end{eqnarray*}
By Lemma \ref{lem:Approx-by-OrbitMeasure} and
(\ref{Eq-choose-n_k}),
$$d\big(\Upsilon_{n_kN_k}(T^{\sum_{j=1}^{k-1} (n_jN_j+t_jl_jL_j)}y),\alpha_{k+1}\big)\leq \zeta_k+2\epsilon_k+d(\alpha_k,\alpha_{k+1})$$
Thus, by
(\ref{diameter-of-Borel-pro-meas}), (\ref{Eq-choose-N_k2222222}) and (\ref{Eq-choose-n_k}),
\begin{eqnarray*}
& &d(\Upsilon_{M_l}(y),\alpha_{k+1}) \\
&\leq&  \frac{\sum_{j=1}^{k-1} (n_jN_j+t_jl_jL_j)}{M_l} d(\Upsilon_{\sum_{j=1}^{k-1} (n_jN_j+t_jl_jL_j)}(y),\alpha'_{M_l})\\
& & + \frac{  n_kN_k }{M_l} d(\Upsilon_{n_kN_k}(T^{\sum_{j=1}^{k-1} (n_jN_j+t_jl_jL_j)}y),\alpha_{k+1}) + \frac{  t_kl_kL_k }{M_l}\times 1\\
& &+ d\big(\Upsilon_{  M_l- \sum_{j=1}^{k} (n_jN_j+t_jl_jL_j)}(T^{   \sum_{j=1}^{k} (n_jN_j+t_jl_jL_j)}y),\alpha_{k+1}\big)\\
&\leq&  \frac{\sum_{j=1}^{k-1} (n_jN_j+t_jl_jL_j)}{\sum_{j=1}^{k} (n_jN_j+t_jl_jL_j) } \times 1  + 1 \times(\zeta_k+2\epsilon_k+d(\alpha_k,\alpha_{k+1}))+ \frac{   t_kl_kL_k }{n_k}
\\ & &
+  2\zeta_{k+1}+2\epsilon_{k+1} \\
&\leq & \zeta_k+\zeta_k+2\epsilon_k+d(\alpha_k,\alpha_{k+1})+\zeta_k+2\zeta_{k+1}+2\epsilon_{k+1}.
\end{eqnarray*}
Since $\zeta_k,\epsilon_k,d(\alpha_k,\alpha_{k+1}) $ all converges to zero as $k$ goes to zero, this proves item (3).

 (4)   As said in the proof  of \cite[Lemma 5.1, item 4]{PS} on Page 946,
the details of the construction are unimportant and in fact  Pfister and Sullivan proved that
\begin{Lem}\label{Lem-PS-Lem5.1-item4}
If $\{n_p\}$ is a strictly increasing sequence of natural numbers such that
$\lim_{p\rightarrow \infty}\frac{M_{n_{p+1}}}{M_{n_p}}=1$ and
$\# \Gamma'_{n_p+1}\times \# \Gamma'_{n_p+2}\cdots \times\# \Gamma'_{n_{p+1}}\geq 2^{h^*(M_{n_{p+1}}-M_{n_p})},$
then $h_{top}(T,G)\geq h^*$.

\end{Lem}

 For $k\geq 1,\, i=0, 1,2,\cdots, N_k-1$, 
 let $$n_{N_1+\cdots+N_{k-1}+i}:=N_1+\cdots+N_{k-1}+t_1+\cdots+t_{k-1}+i$$
Then for any $p\geq 1,$ there is some $k$ so that  $N_1+\cdots+N_{k-1}+t_1+\cdots+t_{k-1}\leq n_p\leq N_1+\cdots+N_{k-1}+t_1+\cdots+t_{k-1}+N_k-1$,  by (\ref{Eq-choose-N_k})
 $$
 1\leq \frac{M_{n_{p+1}}}{M_{n_p}} \leq \frac{M_{n_{p}}+\max\{n_k,\,n_k+t_kl_kL_k\}}{M_{n_p}}=1+\frac{n_k+t_kl_kL_k}{M_{n_p}}
 $$
 $$\leq 1+\frac{n_k+t_kl_kL_k} {\sum_{j=1}^{k-1} (n_jN_j+t_jl_jL_j)}\leq 1+\zeta_k.
 $$
 By (\ref{Eq-choose-n_k2222}) $$\# \Gamma'_{n_p+1}\times \# \Gamma'_{n_p+2}\cdots\times\# \Gamma'_{n_{p+1}-1} \times\# \Gamma'_{n_{p+1}}$$$$ = \# \Gamma_k  \geq 2^{H^*n_k }
 \geq 2^{h^*(n_k+t_kl_kL_k)}\geq 2^{h^*(M_{n_{p+1}}-M_{n_p})}.$$
By Lemma \ref{Lem-PS-Lem5.1-item4} we finish the proof of item (4).

\medskip



 (5) Fix $x\in G.$ By construction, for any fixed $k\geq 1,$ there is $a=a_k$ such that  for any $j=1,\cdots,t_k,$ there is $\Lambda^j\subseteq \Lambda_{l_kL_k}$ $$\max\{d(T^{a+l+{(j-1)l_kL_k}}x,T^{l}x^k_j)|\,l\in \Lambda^j\}\leq\epsilon_k.$$
By (\ref{Eq-minmialrole-Mainproof2222222222222})
$$\frac{\# \Lambda^j}{l_kL_k} \geq 1- \frac{g(l_kL_k)}{l_kL_k}\geq 1-\frac1{4L_k}.$$ Together with  (\ref{Eq-minmialrole-Mainproof}) we get that for  any $j=1,\cdots,t_k$ there is  $p_j\in [0,l_kL_k-1]$  such that $$d(T^{a+p_j+{(j-1)l_kL_k}}x,T^{p_j}x^k_j)\leq \epsilon_k\,\text{ and }\,d( x^k_j,T^{p_j}x^k_j)\leq \epsilon_k.$$  This implies $d(T^{a+p_j+{(j-1)l_kL_k}}x, x^k_j)\leq 2 \epsilon_k$ so that the orbit of $x$ is $3\epsilon_k-$dense in $X$. By arbitrariness of $k,$ one has $x\in Tran$. \qed
\bigskip


If we only have  $g$-almost product property   but do not know the property of uniform separation, then one still has following characterization.

 \begin{Thm}\label{Thm-transtively-convex-saturated}
  Let $T:X\rightarrow X$ be  a continuous map of a compact metric space  $X$ with  $g$-almost product property  property.
   Suppose that there is an invariant measure with full support.     Then $T$ is locally-transitively-convex-saturated and in particular locally-transitively-single-saturated.

\end{Thm}
{\bf Proof.} Similar as the construction  of transitive points in Theorem \ref{Thm-transtively-saturated}, one can adapt the proof of \cite[Theorem 1.5]{DongTian} (or \cite[Theorem 1.2]{PS}) to complete the proof for which ergodic decomposition theorem   replaces the role of uniform separation. Here we omit the details. \qed

\subsection{Locally saturated property}

For any compact connected $K\subseteq M(T,X)$,  let  $G_{K}=\{x\in X|\, M_x=K\},
\,\,\,G^T_{K}
=\{x\in Tran|\, M_x=K
  \}.$  They are saturated set of $K$ and transitively-saturated set of $K$ respectively. 

  \begin{Def}\label{Def-saturated}
 We say that the system $T$ has {\it saturated} property or $T$ is {\it saturated} (simply, S), if  for any  compact connected nonempty set $K \subseteq M (T,X ),$
\begin{eqnarray} \label{eq-saturated-definition}
h_{top} (T,G_K)=\inf\{h_\mu (T)\,|\,\mu\in K\}.
\end{eqnarray}
   We say that the system $T$ has {\it locally-saturated} property or $T$ is {\it locally-saturated} (simply, LS), if  for any  compact connected nonempty set $K \subseteq M (T,X ),$ any nonempty open set $U\subseteq X,$
\begin{eqnarray} \label{eq-locally-saturated-definition}
h_{top} (T,G_K\cap U)=\inf\{h_\mu (T)\,|\,\mu\in K\}.
\end{eqnarray}
    In parallel, one can define    {\it locally-transitively-saturated}   just replacing $ G_K$ by $ G^T_K $ in    (\ref{eq-locally-saturated-definition}). On the other hand,  one can define single-saturated, locally-single-saturated,  locally-transitively-single-saturated   when $K$ is a singleton and define convex-saturated, locally-convex-saturated,   locally-transitively-convex-saturated  respectively  when  $K$   consists of convex sum of two invariant measures.
    \end{Def}

Note that  locally-saturated is stronger than saturated (just taking $U=X$),  saturated is stronger than convex-saturated and the later is stronger than single-saturated.
We will give a basic discussion on the relations of these concepts in Lemma \ref{Lemma-locallytranstive-sameas-transitivesaturated} (see below).
For convenience, we use  $Sat, Loc-Sat, Tran-Sat, Loc-Tran-Sat$  to denote  saturated, locally-saturated,  transitively-saturated, locally-transitively-saturated;
 $ Conv-Sat, Loc-Conv-Sat, Tran-Conv-Sat, Loc-Tran-Conv-Sat $  to denote convex-saturated, locally-convex-saturated,  transitively-convex-saturated, locally-transitively-convex-saturated; and  $Sing-Sat, Loc-Sing-Sat, Tran-Sing-Sat, Loc-Tran-Sing-Sat $  to denote single-saturated, locally-single-saturated,  transitively-single-saturated, locally-transitively-single-saturated, respectively.

 \begin{Lem}\label{Lemma-locallytranstive-sameas-transitivesaturated}
  Let $T:X\rightarrow X$ be  a continuous map of a compact metric space  $X$.
     Then above saturated properties have a following relation:
 $$
\begin{array}{ccccccc}
\text{Tran-Sat}& \Leftrightarrow & \text{Loc-Tran-Sat}& \Rightarrow & \text{Loc-Sat}   & \Rightarrow & \text{Sat}  \\
\Downarrow & & \Downarrow & & \Downarrow& &\Downarrow \\
  \text{Tran-Conv-Sat}	& \Leftrightarrow & \text{Loc-Tran-Conv-Sat}& \Rightarrow & \text{Loc-Conv-Sat} & \Rightarrow & \text{Conv-Sat}  \\
   \Downarrow & & \Downarrow & & \Downarrow& &\Downarrow \\
   \text{Tran-Sing-Sat}	& \Leftrightarrow & \text{Loc-Tran-Sing-Sat}& \Rightarrow & \text{Loc-Sing-Sat} & \Rightarrow & \text{Sing-Sat}  \\
    \end{array}
$$



\end{Lem}

Before proving we need a following basic result, which is also useful for the proof of our main theorems.

\begin{Lem}\label{Lemma-locallytranstive-sameas-transitive2017}
  Let $T:X\rightarrow X$ be  a continuous map of a compact metric space  $X$, $B\subseteq X$ be invariant and $U\subseteq X$ be a nonempty open set.
     Then $h_{top} (T,B\cap Tran\cap U)=h_{top}(T, B\cap Tran). $

\end{Lem}

{\bf Proof.}  By (\ref{Y1containY2}), the part $'\leq'$ is obvious. Now we start to consider the part $'\geq'.$

 Let $B^T=B\cap Tran.$   Notice that by  (\ref{Bowen-fY=Y}) and (\ref{Y1containY2})
  for any $n\geq 1,$ $$ h_{top} (T,f^{-n} (U\cap B^T))=h_{top} (T,f^nf^{-n} (U\cap B^T))\leq h_{top} (T,  U\cap B^T)$$ and
    by the definition of transitivity and invariance of $B$ and $Tran$
 $$B^T=B^T\cap( \bigcup_{n\geq 0}\,f^{-n} U ) \subseteq \bigcup_{n\geq 0}\,f^{-n} (U\cap B^T).$$ Thus  by (\ref{Bowen-disjointSet-entropy})
 $$h_{top}(T,B^T)\leq h_{top}(T,\bigcup_{n\geq 0}\,f^{-n} (U\cap B^T))$$
 $$=\sup_{n\geq 0}h_{top}(T, f^{-n} (U\cap B^T))\leq  h_{top} (T,  U\cap B^T).$$

\begin{Rem} We say $\{Z_i\}$ has {\it locally full entropy gaps} with respect to $Y\subseteq X$ (simply, LFEG w.r.t $Y$) if for any nonempty open set $U\subseteq X,$
$$ h_{top} (T, (Z_{i+1}\setminus Z_i)\cap Y\cap U )=h_{top}(T, Y) \,\,\,\text{ for all } 1\leq i<k.$$ Since the results in this paper are all restricted on transitive points so that by  Lemma \ref{Lemma-locallytranstive-sameas-transitive2017} all results on full entropy gaps can be stated as locally full entropy gaps.

\end{Rem}

{\bf Proof of Lemma \ref{Lemma-locallytranstive-sameas-transitivesaturated}.}  Given a compact and connected subset $K\subseteq M(X,T,) $ we only need to prove that \\
(1)
 $h_{top} (T,G^T_K)=\inf\{h_\mu (T)\,|\,\mu\in K\}\Leftrightarrow h_{top} (T,G^T_K\cap U)
 =\inf\{h_\mu (T)\,|\,\mu\in K\} \text{ holds for any  nonempty open set } U\subseteq X.$\\
 (2) $h_{top} (T,G^T_K\cap U)
 =\inf\{h_\mu (T)\,|\,\mu\in K\} \text{ holds for any  nonempty open set } U\subseteq X $ \\
 $\Rightarrow h_{top} (T,G_K\cap U)
 =\inf\{h_\mu (T)\,|\,\mu\in K\} \text{ holds for any  nonempty open set } U\subseteq X $\\
  $ \Rightarrow h_{top} (T,G_K)
 =\inf\{h_\mu (T)\,|\,\mu\in K\}.$

For part (1), letting $B=G_K$, then it is invariant and by Lemma \ref{Lemma-locallytranstive-sameas-transitive2017}  we know that $h_{top} (T,G^T_K)=  h_{top} (T,G^T_K\cap U).$ This implies part (1).


For part (2), on one hand, if $  h_{top} (T,G^T_K\cap U)
  =\inf\{h_\mu (T)\,|\,\mu\in K\}$   holds for any  nonempty open set  $ U\subseteq X,$ then    $h_{top} (T,G_K\cap U)\geq h_{top} (T,G^T_K\cap U)=\inf\{h_\mu (T)\,|\,\mu\in K\}$. Recall that in general  $h_{top} (f,G_K)\leq \inf\{h_\mu (f)\,|\,\mu\in K\}$ (see  \cite[Theorem 4.1 (3)]{PS2005}), thus $h_{top} (T,G_K\cap U)=\inf\{h_\mu (T)\,|\,\mu\in K\}$.   On the other hand, if $  h_{top} (T,G_K\cap U)
  =\inf\{h_\mu (T)\,|\,\mu\in K\}$   holds for any  nonempty open set  $ U\subseteq X,$ then taking $U=X,$ we have $h_{top} (T,G_K)=h_{top} (T,G^T_K\cap U)=\inf\{h_\mu (T)\,|\,\mu\in K\}$.
\qed

\begin{Rem} If $G^T_K\neq \emptyset,$ then there is transitive point so that the system is transitive and thus is surjective. So  for $n\geq 1,$
$$ h_{top} (T,f^{-n} (U\cap G^T_K))=h_{top} (T,f^nf^{-n} (U\cap G^T_K))= h_{top} (T,  U\cap G^T_K).$$  By definition and   the system being surjective, one also see that  $f^{\pm 1}G_K=  G_K,$ $ f^{\pm 1}Tran =  Tran, $ $f^{\pm 1}G^T_K=  G^T_K$.

\end{Rem}




\section{Multifractal Analysis, Transitive and Banach Recurrence}\label{section-prove-banach-multifractal}

\subsection{Transitive and Banach  Recurrence}\label{section-20150108-recurrence}


   A point $x$ is called {\it quasi-generic} for some measure $\mu,$ if there are two sequences of positive integers $\{a_k\}$ and $\{b_k\}$
  with $b_k\geq a_k,\,\,\lim_{k\rightarrow\infty}b_k-a_k=\infty$ such that $$\lim_{k\rightarrow\infty}\frac{1}{b_k-a_k+1}\sum_{j=a_k}^{b_k}\delta_{T^j(x)}=\mu$$
   in weak$^*$ topology.  This concept is from \cite[Page65]{Furst}.  Let $M_x^*=\{\mu\in M(T,X): \,x \text{ is quasi-generic for } \mu\}$
   and  $C^*_x=\overline{\cup_{m\in M^*_x  }S_m}$. Let $C_x=\overline{\cup_{m\in M_x  }S_m}$. Note that $M_x\subseteq M^*_x\subseteq M(T,X)$ and $C_x\subseteq C^*_x\subseteq X.$  From \cite[Proposition 3.9]{Furst}
   it is known that  $M^*_x$ is always nonempty, compact and connected.
 It is not difficult to show that (see    \cite{Zhou-center-measure,Zhou93,Zhou95,ZF,Huang-Wang}):

\begin{Lem}\label{lem-zhou-centerminimal} For $\forall\,\, x\in X,$
 \begin{eqnarray}\label{C_xinOmega}
\,\,\,C_x\subseteq C^*_x\subseteq \omega_T (x);
\end{eqnarray}
\begin{eqnarray}\label{QW}
\,x\in QW  \Leftrightarrow x\in C_x\Leftrightarrow x\in\omega_T (x)=C_x 
\end{eqnarray}
 \begin{eqnarray}\label{Banachrecurrent}
\,x\in BR  \Leftrightarrow x\in C^*_x\Leftrightarrow x\in\omega_T (x)=C^*_x 
\end{eqnarray}

\end{Lem}


\begin{Lem}\label{Lem-Cstar=X}   For $x\in X,$ if $C^*_x=X,$  then  $x\in BR \cap Tran  .$

\end{Lem}

{\bf Proof.}  By (\ref{C_xinOmega}) $\omega_T (x)=X$  and $x\in X=  C^*_x$ so that by (\ref{Banachrecurrent}), $x\in BR\cap Tran  .$  \qed

\begin{Lem}\label{Lem-Tran-in-Cstar}
For any $x\in X,$
$C_x^*=\overline{\cup_{\nu\in M_{erg}(T,\omega_T(x))}S_\nu}= \overline{\cup_{\nu\in M(T,\omega_T(x))}S_\nu}.$ 
Moreover, there is an invariant measure $\mu\in  M(T,\omega_T(x))$  so that $S_\mu=C^*_x.$ (We emphasize that it is unknown whether there is an invariant measure $\mu\in  M^*_x$  so that $S_\mu=C^*_x$).

\end{Lem}

{\bf Proof.} It is obvious that $ {\cup_{\nu\in M_{erg}(T,\omega_T(x))}S_\nu}\subseteq \overline{\cup_{\nu\in M(T,\omega_T(x))}S_\nu}.$ By ergodic decomposition theorem, we know that for any
invariant measure $\mu,$ $$\mu( \overline{\cup_{\nu\in M_{erg}(T,\omega_T(x))}S_\nu})=1$$ so that $S_\mu\subseteq \overline{\cup_{\nu\in M_{erg}(T,\omega_T(x))}S_\nu}$ and thus  $  {\cup_{\nu\in M(T,\omega_T(x))}S_\nu}\subseteq \overline{\cup_{\nu\in M_{erg}(T,\omega_T(x))}S_\nu}$.

From \cite [Proposition 3.9, Page 65] {Furst} we know that for a point $x_0$ and an ergodic measure $\mu_0\in M(\omega_T(x_0),T)$, $x_0$ is quasi-generic for $\mu_0.$ This implies that for any $x\in X,$ $M_{erg}(T,\omega_T(x))\subseteq M^*_x$. So $ \overline{\cup_{\mu\in M_{erg}(T,\omega_T(x))}S_\mu}\subseteq C^*_x \subseteq \overline{\cup_{\mu\in M(T,\omega_T(x))}S_\mu}.$

Take a sequence of invariant measures $\mu_i\in M(T,\omega_T(x))$ dense in $M(T,\omega_T(x))$
 such that $ \overline{\cup_{i\geq 1}S_{\mu_i}}= \overline{\cup_{\mu\in M(T,\omega_T(x))}S_\mu}.$ Let $\mu=\sum_{n\geq 1}\frac1{2^n}\mu_n.$ Then $\mu(\cup_{i\geq 1}S_{\mu_i})=1$ so that $S_\mu=\overline{\cup_{\nu\in M(T,\omega_T(x))}S_\nu}=C_x^*.$\qed

\medskip

\begin{Lem}\label{Lem-Tran-in-BR000000} For any $x\in Rec,$
 $$x\in BR\Leftrightarrow \text{ there is an invariant measure $\mu\in  M(T,\omega_T(x))$  so that $S_\mu=\omega_T(x).$} $$
\end{Lem}

{\bf Proof.} The part $\Rightarrow$: By (\ref{Banachrecurrent})  $x\in  \omega_T (x)=C^*_x$. By Lemma \ref{Lem-Tran-in-Cstar}, there is an invariant measure $\mu\in  M(T,\omega_T(x))$  so that $S_\mu=C^*_x.$ So $S_\mu=C^*_x=\omega_T(x).$


 The part $\Leftarrow$: If there is an invariant measure with $S_\mu=\omega_T(x)\ni x$, then  by (\ref{C_xinOmega}) $ C_x^* \subseteq \omega_T(x) $ and by Lemma \ref{Lem-Tran-in-Cstar},
 $C_x^*=  \overline{\cup_{\nu\in M(T,\omega_T(x))}S_\nu} \supseteq S_\mu=\omega_T(x)$.  
 By (\ref{Banachrecurrent}) $x\in BR$. \qed


By Lemma \ref{Lem-Tran-in-BR000000}, if $x\in Tran,$ then $$x\in   BR \Leftrightarrow \text{ there is  $\mu\in  M(T,X)$  so that $S_\mu=X.$} $$ 
Thus, we have
\begin{Lem}\label{Lem-Tran-in-BR} Suppose $(X,T)$ is transitive. Then  
 $Tran\subseteq BR \Leftrightarrow Tran\cap BR\neq \emptyset \Leftrightarrow$ there is an invariant measure with full support.

\end{Lem}






\medskip

{  Let $BV:=\{x\in BR|\,\exists\,\,\mu\in M^*_x\,\,s.t.\,\, S_\mu=C^*_x\}.$
\begin{Lem}\label{Lem-Tran-in-BanV} If ergodic measures are dense in the space of invariant measures (or $T$ has entropy-dense property), then  for any $x\in Tran,$ $M^*_x=M(T,X)$. If further $T$ has an invariant measure with full support, then  $Tran\subseteq BV$. 

\end{Lem}

{\bf Proof.} From \cite [Proposition 3.9, Page 65] {Furst} we know that for a point $x_0$ and an ergodic measure $\mu_0\in M(\omega_T(x_0),T)$, $x_0$ is quasi-generic for $\mu_0.$ This implies that for any $x\in Tran,$ $M_{erg}(T,X)\subseteq M^*_x$. By assumption of density of ergodic measures, $M^*_x=M(T,X)$. If further $T$ has an invariant measure $\mu$ with full support, then $\mu\in M_x^*$
and $C^*_x=S_\mu=X=\omega_T(x)$ and by (\ref{Banachrecurrent}) $x\in BV.$ \qed}


For a set $K\subseteq M(T,X),$ define $C_K=\overline{\cup_{\omega\in K}S_\omega}.$ Recall the notions that $G_{K}=\{x\in X|\, M_x=K\},
\,\,\,G^T_{K}
=\{x\in Tran|\, M_x=K
  \}.$ Let $$    V^* =\{ \,x\in Rec  |\,\exists \,\mu\in M_x  \text{ such that } S_\mu=C_x\},\,$$
  $$ W^* =\{\, x\in  Rec \,|\,S_\mu=C_x  \text{ for every }\mu \in M_x \}.$$
\begin{Lem}\label{lem-various-GK-2016} Let $K$ be a compact connected subset of $M(T,X)$. \\
(1) If for any $\omega\in K,$ $S_\omega=C_K ,$ then
$G_K \subseteq W^*.$\\
(2) If there are two measures $\omega_i\in K\,(i=1,2),$ $S_{\omega_1}\subsetneq S_{\omega_2}=C_K ,$ then
$G_K \subseteq (V^*\setminus W^*).$\\
(3) If  for any $\omega\in K,$ $S_\omega\neq C_K, $ then
$ G_K\subseteq  X\setminus V^*.$\\
(4)
If  $T|_{C_K}$ is not a transitive subsystem,  then
$G_K\subseteq X\setminus QW.$ \\
(5) If $C_K\neq X$,  then $G^T_K\subseteq   Tran\setminus QW.$

\end{Lem}

{\bf Proof.} 
  The proofs of (1)-(3) are not difficult. 
     We only give the proof  of (4) and (5).

 For case (4), by contradiction there is $x\in G_K\cap QW$, then $C_x=C_K$ and by (\ref{QW}) $x\in\omega_T(x)= C_x$ so that $x\in\omega_T(x)= C_K.$ It means $T|_{C_K}$ is transitive, it contradicts the assumption.

  For case (5), by contradiction there is $x\in G^T_K\cap QW$, then $x\in G^T_K$ implies that $C_x=C_K\neq X,\,\omega_T(x)=X$ so that $C_x\neq \omega_T(x)$ and by (\ref{QW}) $x\in   QW$ implies that $x\in\omega_T(x)= C_x$. It is a contradiction.
  \qed

\subsection{Proof of Theorem \ref{Thm--RoughResult-multifrcalanalysis}: Irregular set and Level sets}
Recall that the system $T$ has {\it single-saturated} property or $T$ is {\it single-saturated}, if  for any    $\mu\in M (T,X ),$
\begin{eqnarray} \label{eq-single-saturated-definition}
h_{top} (T,G_\mu)= h_\mu .
\end{eqnarray}
Now let us state a result on variational principle of $R_\alpha(a) $. Let $$H_\alpha(a)=\sup \{h_{\rho}  |\,\,\rho\in M (T,X)\,\,and\,\,\alpha (\rho)=a\} .$$

\begin{Prop}\label{Main-prop-levelsets-0} 
   Suppose that $T$ is transitively-single-saturated.
 Let $\alpha:M(T,X)\rightarrow \mathbb{R}$ be a continuous function.     Then for any real number $a\in  L_\Phi ,$
$$ h_{top} (T,R_\alpha(a))=h_{top} (T,R_\alpha(a)\cap Tran  )=H_\alpha(a).$$
If further $T$ has positive topological entropy and $Int(L_\Phi)\neq \emptyset,$ then  for any real number $a\in  Int(L_\Phi) ,$
 $$  h_{top} (T,R_\alpha(a))=h_{top} (T,R_\alpha(a)\cap Tran  )=H_\alpha(a)>0.$$

\end{Prop}
{\bf Proof.}  For any real number $t \geq 0$, define the (maybe empty) set
 $$Q(t) := \{x:\exists \mu\in M_x \,\,s.t.\,\,h_\mu(f)\leq t\}.$$
 From \cite[Theorem 2]{Bowen1}:
   $h_{top}(f, Q(t))\leq t$.
Let  $t= H_\alpha(a)$, notice that  $R_\alpha(a) \subseteq Q(t)$ and thus $h_{top} (T,R_\alpha(a)\cap Tran  )\leq h_{top}(T,R_\alpha(a))\leq h_{top}(T, Q(t))\leq t$.

On the other hand, for any invariant measure $\rho$  with $\alpha(\rho)=a,$   note that $ G^T_\rho\subseteq R_\alpha(a)\cap Tran.$ By assumption of transitively-single-saturated,
 $h_{top}(T,R_\alpha(a)\cap Tran)\geq  h_{top}(T,G^T_\rho)=h_\rho.$

 Suppose that further $T$ has positive topological entropy and $Int(L_\Phi)\neq \emptyset,$ fix $a\in  Int(L_\Phi) .$
  By classical variational principle, we can take an ergodic measure $\rho_1$ such that $h_{\rho_1}>0$. If $\alpha(\rho_1)=a,$ then $H_\alpha(a)>0.$  If $\alpha(\rho_1)\neq a,$ without loss of generality, we may assume that $\alpha(\rho_1)< a$. Since $a\in  Int(L_\Phi) ,$ we can take another invariant measure $\rho_2$ such that  $\alpha(\rho_2)> a$. Then one can take suitable $\theta\in(0,1)$ such that $\rho:=\theta\rho_1+(1-\theta)\rho_2$ satisfies that  $\alpha(\rho)=a.$ Note that $h_\rho\geq \theta h_{\rho_1}>0,$  then $H_\alpha(a)>0.$
  \qed

\medskip






From Birkhoff ergodic theorem,  $QR $ has totally full measure and note that $QR\cap I_\alpha=\emptyset$. Thus 
 $I_\alpha $ has zero measure for any invariant measure. 
However, in recent several
   years many people have focused on studying the dynamical complexity of  irregular set from different sights, for example,
    in the sense of dimension theory and topological entropy (or pressure) etc.  Pesin and Pitskel   \cite{Pesin-Pitskel1984} are the first to notice the phenomenon of the irregular set carrying full topological entropy in the case of the full shift on two symbols. Barreira, Schmeling, etc. studied   the irregular set  in the setting
of subshifts of finite type and beyond, see    \cite{Barreira-Schmeling2000,Bar,TV} etc.
     Recently,  Thompson shows in    \cite{To2010,Tho2012} that  every $\alpha-$irregular set  $I_\alpha  $ (in this case $\alpha:=\int\Phi d\mu$ for additive functions of $\Phi:X\rightarrow \mathbb{R}$) either is empty or carries full topological entropy
 (or pressure) when the system satisfies (almost) specification, which is inspired from    \cite{PS} by Pfister
     and Sullivan and    \cite{TV} by Takens and Verbitskiy. Now we sate a result for a continuous function $\alpha: M(T,X)\rightarrow \mathbb{R}$.






\begin{Lem}\label{Lem-IC-notempty} Suppose that $T$ is transitively-convex-saturated and has positive topological entropy.  Let $\alpha:M(T,X)\rightarrow \mathbb{R}$ be a continuous function.
     If      $$\inf_{\mu\in M(T,X)}\alpha( \mu)<\sup_{\mu\in  M(T,X)}\alpha( \mu),$$  then  $h_{top}(T,I_\alpha)\geq h_{top}(T,I_\alpha\cap Tran)>0$.  
If further  $\alpha:M(T,X)\rightarrow \mathbb{R}$   satisfies
$[A.3]$,
  then  $h_{top}(T,I_\alpha)= h_{top}(T,I_\alpha\cap Tran)=h_{top}$.
\end{Lem}

 {\bf Proof.}    Fix $\epsilon\in(0,h_{top}).$ By classical variational principle, we can take an ergodic measure $\nu$ such that $h_{\nu}>h_{top}-\epsilon>0$.  By assumption we can another invariant measure $ \nu_1$ such that  $\alpha(\nu_1)\neq \alpha(\nu)$.
   Then by continuity of $\alpha$ there is  $\theta\in(0,1)$ such that $\rho:=\theta\nu+(1-\theta)\nu_1$ satisfies that  $\alpha(\rho)\neq \alpha(\nu).$ Remark that $h_\rho\geq \theta h_\nu>0$. Let $K=\{t\nu+(1-t)\mu|\,t\in[0,1]\}.$ Since $T$ is transitively-convex-saturated, one can get $ h_{top}(T,G^T_K)=\min\{h_\nu,h_\rho\}>0$. Note that $G^T_K\subseteq I_\alpha\cap Tran$. Thus $h_{top}(T,I_\alpha\cap Tran)\geq h_{top}(T,G^T_K) >0$.

  If further  $\alpha:M(T,X)\rightarrow \mathbb{R}$   satisfies
$[A.3]$,  then above $\theta\in(0,1)$ can be chosen very close to 1 such that $\rho:=\theta\nu+(1-\theta)\nu_1$ satisfies that  $h_\rho\geq \theta h_\nu>h_{top}-\epsilon$. Then one has $h_{top}(T,I_\alpha\cap Tran)\geq h_{top}(T,G^T_K)=\min\{h_\nu,h_\rho\}>h_{top}-\epsilon$.
 \qed

{\bf Proof of Theorem \ref{Thm--RoughResult-multifrcalanalysis}. } It follows from  Lemma \ref{Lem-IC-notempty} and Proposiition \ref{Main-prop-levelsets-0}. \qed




\begin{Prop}\label{Prop-Fact-sup-H-alpha-a}    Let $\alpha:M(T,X)\rightarrow \mathbb{R}$ be a continuous function.  Then $h_{top}(T) =h_{top}(R_\alpha).$
If  further   $Int (L_\alpha) \neq \emptyset$ and $\alpha$ satisfies $[A.3],$ then
  $$h_{top}(T)=\sup_{a\in Int (L_\alpha)} H_\alpha(a)=h_{top}(R_\alpha).$$

  \end{Prop}

{\bf Proof.} Note that  $QR\subseteq R_\alpha$ so that $\mu(R_\alpha)=1$ for any invariant measure $\mu.$ By \cite[Theorem 3]{Bowen1} $h_{top} (T,\Gamma)\geq h_\mu(f)$ for any $\Gamma$ with $\mu(\Gamma)=1$. Thus
$h_{top}(R_\alpha)=\sup_{\mu\in M(T,X)} h_\mu=h_{top}(T)  \geq \sup_{a\in Int (L_\alpha)} H_\alpha(a)$. Now we only need to prove  $h_{top}(T)\leq \sup_{a\in Int (L_\alpha)} H_\alpha(a)$.

By  classical variational principle, there is an ergodic measure $\mu$ such that $h_\mu>h_{top}-\epsilon.$ If $\alpha(\mu)\in Int (L_\alpha),$ take $\omega=\mu.$ Otherwise, take invariant measure $\nu$ such that $\alpha(\nu)\neq \alpha(\mu)$ and $\alpha(\nu)\in Int (L_\alpha).$ By condition $[A.3]$ and continuity of $\alpha,$ one can choose $\theta\in(0,1)$ close to 1 such that $\omega=\theta\mu+(1-\theta)\nu$ satisfies $\alpha(\omega)\in Int (L_\alpha)$ and  $h_\omega\geq \theta h_\mu>h_{top}-\epsilon.$ Thus 
 $ \sup_{a\in Int (L_\alpha)} H_\alpha(a)\geq H_\alpha(\alpha(\omega)) >h_{top}-\epsilon.$ \qed




\subsection{Proof  of Theorem \ref{Maintheorem111-2017-irregular}: Gap-sets w.r.t. $I_\alpha$} \label{proof-main-irregular}

\begin{mainlemma}\label{Main-Lemma-irregular}

 Suppose that  $T$ has entropy-dense property. Let $\alpha:M(T,X)\rightarrow \mathbb{R}$ be a continuous function satisfying $[A.3]$ and   $Int (L_\alpha) \neq \emptyset$. Then for any
 $\epsilon>0,$ any integer $k\geq 2,$ there exist ergodic measures $\rho_1,\rho_2,\cdots,\rho_k,$ such that
 \begin{eqnarray*}
&(1)&  h_{\rho_i}>h_{top}-\epsilon,\,i=1,2,\cdots,k,\\
 &(2)& S_{\rho_i}\cap S_{\rho_j}\neq \emptyset,\,1\leq i<j\leq k,\\
 &(3)& \alpha(\rho_1)<\alpha(\rho_2)<\cdots<\alpha(\rho_k).
 \end{eqnarray*}

\end{mainlemma}

{\bf Proof.} Fix $\epsilon>0.$ By classical variational principle, we can take an ergodic measure $\nu$ such that $h_{\nu}>h_{top}-\frac12\epsilon$.  By assumption we can choose another invariant measure $ \nu_1$ such that  $\alpha(\nu_1)\neq \alpha(\nu)$. Then there is   $\theta_0\in(0,1)$   very close to 1  such that $h_{\theta\nu+(1-\theta)\nu_1}\geq \theta h_\nu\geq   \theta_0 h_\nu>h_{top}-\frac12\epsilon$ hold for any $\theta\in(\theta_0,1).$
   By condition $[A.3]$ and continuity of $\alpha$ there are  $\theta_0<\theta_1<\theta_2<\cdots<\theta_k<1$    such that  $\tau_i:=\theta_i\nu+(1-\theta_i)\nu_1$ satisfies that  $h_{\tau_i}\geq \theta_i h_\nu>h_{top}-\frac12\epsilon,\,i=1,2,\cdots,k$ and  $\alpha(\tau_1)<\alpha(\tau_2)<\cdots<\alpha(\tau_k)$ or $\alpha(\tau_1)>\alpha(\tau_2)>\cdots>\alpha(\tau_k).$ We may assume that $\alpha(\tau_1)<\alpha(\tau_2)<\cdots<\alpha(\tau_k)$ (Otherwise, by reverse order to change).

By continuity of $\alpha$, take  neighborhood $G_i$ of $\tau_i$ in $M(T,X)$ ($i=1,2,\cdots,k$) such that $G_i\cap G_j=\emptyset,\,\sup_{\tau\in G_i}\alpha(\tau)<\inf_{\tau\in G_j}\alpha(\tau), \,1\leq i<j\leq k.$  By  entropy-dense property,    there exist  $\rho_i\in M_{erg}(T,  X)$ ($i=1,2,\cdots,k$) such that $h_{\rho_i}>h_{\tau_i}-\frac12\epsilon$ and $M(T,  S_{\rho_i})\subseteq G_i$. Then these ergodic measures $\rho_1,\rho_2,\cdots,\rho_k,$ are required. \qed

\bigskip

  Let $ Tran^\# :=\{x\in Tran |\, C_x\subsetneq C^*_x\}$,
 \begin{eqnarray*}   W_T^\# &:=&\{\, x\in  Tran^\# \,|\,S_\mu=C_x  \text{ for every }\mu \in M_x \},
   \\   V_T^\#  &:=&\{ \,x\in  Tran^\#  |\,\exists \,\mu\in M_x  \text{ such that } S_\mu=C_x\},\\
   S  &:=&\{ \,x\in  X  |\,\cap_{\mu\in M_x} S_\mu\neq \emptyset\};\\
   Tran_1  &:=&  W_T^\#,
   Tran_2   :=  V_T^\#\cap S ,\,\,\,\, Tran_3   :=  V_T^\#, \\
   Tran_4  &:=& V_T^\#\cup (Tran^\#\cap S) ,\,\,\,\,\, Tran_5   :=  Tran^\#,
          \end{eqnarray*}

 \begin{Thm}\label{mainthm--Maintheorem111-2017-irregular}

 Suppose that $(X,\,T)$   has transitively-convex-saturated property and  entropy-dense property. Let $\alpha:M(T,X)\rightarrow \mathbb{R}$ be a continuous function  satisfying   $[A.3]$ and $Int(L_\alpha)\neq \emptyset.$
   Then   $\{\emptyset,  Tran_1,Tran_2,Tran_3,Tran_4,Tran_5 \}$ has   full entropy gap w.r.t. $I_\alpha $.

 \end{Thm}

{\bf Proof.} Fix $\epsilon>0.$ Take $\rho_1,\rho_2,\rho_3,\rho_4$ satisfying Lemma \ref{Main-Lemma-irregular}. By condition $[A.3]$ of $\alpha,$ one can take suitable $\theta_i\in (0,1)$ ($i=1,2$)
such that $\nu_i=\theta_i\rho_1+(1-\theta_i)\rho_{2}$ satisfies that $\alpha(\nu_1)\neq \alpha(\nu_2)$. Note that $S_{\nu_i}=S_{\rho_1}\cup S_{\rho_2}$ and  $h_{\nu_i}=\theta_ih_{\rho_1}+(1-\theta_i)h_{\rho_{2}}>h_{top}-\epsilon.$ Take   suitable $\theta \in (0,1)$
such that $\omega=\theta \rho_1+(1-\theta )\rho_{3}$  satisfies that $\alpha(\omega)\neq \alpha(\rho_1)$. Note that $S_{\omega}=S_{\rho_1}\cup S_{\rho_3}$ and  $h_{\omega}=\theta h_{\rho_1}+(1-\theta )h_{\rho_{3}}>h_{top}-\epsilon.$
 Let \begin{eqnarray*}
&K_1:=& \{t\nu_1+(1-t)\nu_2|\,t\in[0,1]\},\\
&K_2:=& \{t\rho_1+(1-t)\omega|\,t\in[0,1]\},\\
&K_3:=& \{t\rho_1+(1-t)\rho_2|\,t\in[0,1]\},\\
&K_4:=& \{t\rho_1+(1-t)\nu_1|\,t\in[0,1]\}\cup \{t\rho_1+(1-t)\omega|\,t\in[0,1]\},\\
&K_5:=& \{t\rho_1+(1-t)\rho_2|\,t\in[0,1]\}\cup \{t\rho_1+(1-t)\rho_3|\,t\in[0,1]\}.
 \end{eqnarray*}
By  transitively-convex-saturated property, $$h_{top}(G^T_{K_i})=\inf_{\mu\in K_i}h_\mu\geq \min\{h_{\rho_1},h_{\rho_2},h_{\rho_3},h_{\nu_1},h_{\nu_2},h_{\omega}\} >h_{top}-\epsilon.$$
Note that $K_i\subseteq I_\alpha$ and $C_{K_i}\subseteq \cup_{i=1}^3S_\mu\subsetneq \cup_{i=1}^4S_\mu\subseteq X$. By Lemma \ref{lem-various-GK-2016} $G^T_{K_i}\subseteq   Tran\setminus QW$ and
$G^T_{K_i}\subseteq Tran_i\setminus Tran_{i-1}$ where $Tran_0=\emptyset,\,i=1,2,\cdots,5.$
 Then we complete the proof. \qed

{\bf Proof  of Theorem \ref{Maintheorem111-2017-irregular}.} By Lemma \ref{Lem-Tran-in-BR}, $Tran\subseteq BR$ so that $Tran_i\subseteq BR_i.$ Thus one can use Theorem \ref{mainthm--Maintheorem111-2017-irregular} to complete the proof. \qed

\subsection{Proof  of Theorem \ref{Maintheorem222-2017-levelsets}: Gap-sets w.r.t. level set $R_\alpha(a)$} \label{proof-main-levelset}

\begin{mainlemma}\label{Main-Lemma-levelset}

 Suppose that  $T$ has entropy-dense property. Let $\alpha:M(T,X)\rightarrow \mathbb{R}$ be a continuous function satisfying $[A.1]$  and   $Int (L_\alpha) \neq \emptyset$. Then for any
 $\epsilon>0,$ any $ a\in Int(L_\alpha)$, any integers $k\geq 2, \,l\geq 2$ there exist ergodic measures $\rho_1,\rho_2,\cdots,\rho_{k+l},$ such that
 \begin{eqnarray*}
&(1)&  h_{\rho_i}>H_\alpha(a)-\epsilon,\,i=1,2,\cdots,k+l,\\
 &(2)& S_{\rho_i}\cap S_{\rho_j}\neq \emptyset,\,1\leq i<j\leq k+l,\\
 &(3)& \alpha(\rho_1)<\alpha(\rho_2)<\cdots\alpha(\rho_k)<a<\alpha(\rho_{k+1})<\cdots\alpha(\rho_{k+l}).
 \end{eqnarray*}

\end{mainlemma}

{\bf Proof.} Fix $\epsilon>0$ and  $ a\in Int(L_\alpha)$. Take an ergodic measure $\nu$ with $\alpha(\nu)=a$ such that $h_{\nu}>H_\alpha(a)-\frac12\epsilon$.  By assumption we can choose another two invariant measure $ \nu_1,\nu_2$ such that
 $\alpha(\nu_1)< a=\alpha(\nu)< \alpha(\nu_2)$. Then there is   $\theta_0\in(0,1)$   very close to 1  such that $h_{\theta\nu+(1-\theta)\nu_i}\geq \theta h_\nu\geq   \theta_0 h_\nu>H_\alpha(a)-\frac12\epsilon$ hold for $i=1,2$ and any $\theta\in(\theta_0,1).$
   By condition $[A.1]$ and continuity of $\alpha$ there are  $\theta_0<\theta_1<\theta_2<\cdots<\theta_k<1$  and $\theta_0<\theta_{k+1}<\theta_{k+2}<\cdots<\theta_{k+l}<1$   such that  $\tau_i:=\theta_i\nu+(1-\theta_i)\nu_1$ satisfies that  $h_{\tau_i}\geq \theta_i h_\nu>H_\alpha(a)-\frac12\epsilon,\,i=1,2,\cdots,k+l$ and  $\alpha(\tau_1)<\alpha(\tau_2)<\cdots<\alpha(\tau_k)<a<\alpha(\tau_{k+1})<\cdots\alpha(\tau_{k+l}).$

By continuity of $\alpha$, take  neighborhood $G_i$ of $\tau_i$ in $M(T,X)$ ($i=1,2,\cdots,k+l$) such that $G_i\cap G_j=\emptyset,\,\sup_{\tau\in G_i}\alpha(\tau)<\inf_{\tau\in G_j}\alpha(\tau), \,1\leq i<j\leq k+l.$  By  entropy-dense property,    there exist  $\rho_i\in M_{erg}(T,  X)$ ($i=1,2,\cdots,k+l$) such that $h_{\rho_i}>h_{\tau_i}-\frac12\epsilon$ and $M(T,  S_{\rho_i})\subseteq G_i$. Then these ergodic measures $\rho_1,\rho_2,\cdots,\rho_{k+l},$ are required. \qed

\begin{mainlemma}\label{Main-Lemma-levelset222222222}

 Suppose that  $T$ has entropy-dense property. Let $\alpha:M(T,X)\rightarrow \mathbb{R}$ be a continuous function satisfying $[A.1]$  and   $Int (L_\alpha) \neq \emptyset$. Then for any
 $\epsilon>0,$ any $ a\in Int(L_\alpha)$, any integers $k\geq 2, $ there exist invariant  measures $\mu_1,\mu_2,\cdots,\mu_{k },$ such that
 \begin{eqnarray*}
&(1)&  h_{\mu_i}>H_\alpha(a)-\epsilon,\,i=1,2,\cdots,k ,\\
 &(2)& S_{\mu_i}\cap S_{\mu_j}\neq \emptyset,\,1\leq i<j\leq k,\\
 &(3)& \alpha(\mu_i)=a,\,i=1,2,\cdots,k.
 \end{eqnarray*}

\end{mainlemma}

{\bf Proof.} Fix $\epsilon>0$ and  $ a\in Int(L_\alpha)$. Take  ergodic measures $\rho_1,\rho_2,\cdots,\rho_{2k}$ (in the case $l=k$) satisfying
Lemma \ref{Main-Lemma-levelset}. By continuity of $\alpha,$ take suitable $\theta_i\in(0,1)$ such that $\mu_i=\theta_i\rho_i+(1-\theta_i)\rho_{2k+1-i}$ satisfies item (3), $i=1,2,\cdots,k$. Then it is easy to check such $\mu_1,\mu_2,\cdots,\mu_{k }$ are required.\qed

\begin{Thm}\label{mianthm--Maintheorem222-2017-levelsets} 
Suppose that $(X,\,T)$   has transitively-convex-saturated property and  entropy-dense property.
Let $\alpha:M(T,X)\rightarrow \mathbb{R}$ be a continuous function  satisfying  $[A.1]$, $[A.2]$ and $Int(L_\alpha)\neq \emptyset.$ Let $a\in Int(L_\alpha)$.
 Then    $$\{\emptyset, QR\cap Tran_1, Tran_1,Tran_2,Tran_3,Tran_4,Tran_5 \}$$  has   full entropy gap w.r.t. $R_\alpha(a)\cap Tran $.


\end{Thm}

{\bf Proof.} Fix $\epsilon>0.$ Take $\mu_1,\mu_2,\mu_3,\mu_4$ satisfying Lemma \ref{Main-Lemma-levelset222222222}. By condition $[A.2]$ of $\alpha,$ one can take different $\theta_i\in (0,1)$ ($i=1,2$)
such that $\nu_i=\theta_i\mu_1+(1-\theta_i)\mu_{2}$ satisfies that $\nu_1\neq\nu_2,\,\alpha(\nu_1)= \alpha(\nu_2)=a$. Note that $S_{\nu_i}=S_{\mu_1}\cup S_{\mu_2}$ and  $h_{\nu_i}=\theta_ih_{\mu_1}+(1-\theta_i)h_{\mu_{2}}>H_\alpha(a)-\epsilon.$ Take     $\theta \in (0,1)$
such that $\omega=\theta \mu_1+(1-\theta )\mu_{3}$  satisfies that $\alpha(\omega)= \alpha(\mu_1)=a$ by  condition $[A.2]$. Note that $S_{\omega}=S_{\mu_1}\cup S_{\mu_3}$ and  $h_{\omega}=\theta h_{\mu_1}+(1-\theta )h_{\mu_{3}}>H_\alpha(a)-\epsilon.$
 Let \begin{eqnarray*}
&K_0:=& \{\mu_1\},\\
&K_1:=& \{t\nu_1+(1-t)\nu_2|\,t\in[0,1]\},\\
&K_2:=& \{t\mu_1+(1-t)\omega|\,t\in[0,1]\},\\
&K_3:=& \{t\mu_1+(1-t)\mu_2|\,t\in[0,1]\},\\
&K_4:=& \{t\mu_1+(1-t)\nu_1|\,t\in[0,1]\}\cup \{t\mu_1+(1-t)\omega|\,t\in[0,1]\},\\
&K_5:=& \{t\mu_1+(1-t)\mu_2|\,t\in[0,1]\}\cup \{t\mu_1+(1-t)\mu_3|\,t\in[0,1]\}.
 \end{eqnarray*}
By  transitively-convex-saturated property, $$h_{top}(G^T_{K_i})=\inf_{\mu\in K_i}h_\mu\geq \min\{h_{\mu_1},h_{\mu_2},h_{\mu_3},h_{\nu_1},h_{\nu_2},h_{\omega}\} >H_\alpha(a)-\epsilon.$$
Note that $K_i\subseteq R_\alpha(a)$ and $C_{K_i}\subseteq \cup_{i=1}^3S_\mu\subsetneq \cup_{i=1}^4S_\mu\subseteq X$.
 By  Lemma \ref{lem-various-GK-2016}
  $G^T_{K_i}\subseteq   Tran\setminus QW, \,i=0,1,\cdots,5$ and $G^T_{K_0}\subseteq QR\cap Tran_1$, $G^T_{K_1}\subseteq  Tran_1\setminus QR$,
$G^T_{K_i}\subseteq Tran_i\setminus Tran_{i-1}$,  $ i= 2,3,4,5.$
 Then we complete the proof. \qed

{\bf Proof  of Theorem  \ref{Maintheorem222-2017-levelsets}.} By Lemma \ref{Lem-Tran-in-BR}, $Tran\subseteq BR$ so that $  Tran_i\subseteq BR_i, \,i= 1,2,3,4,5.$ Thus one can use Theorem \ref{mianthm--Maintheorem222-2017-levelsets} to complete the proof. \qed

\subsection{Proof of Theorem \ref{Maintheorem000-2017-banachrecurrence} and \ref{Maintheorem333-2017-alpha-regularsets}.}

{\bf Proof of Theorem \ref{Maintheorem000-2017-banachrecurrence}. }
Since the system is not uniquely ergodic, then there exist two different invariant measures $\mu,\nu.$ Then by weak$^*$ topology, there is a continuous function $\phi:X\rightarrow \mathbb{R}$ such that $\int \phi d\mu \neq \int \phi d\nu.$ Define $\alpha:M(T,X)\rightarrow \mathbb{R},\,\tau\mapsto \int \phi d\tau.$ By convex-saturated property, there is $x\in X$ such that $M_x=\{t\mu+(1-t)\nu|\,t\in[0,1]\}.$ Then $x\in I_\alpha$ so that one can use
Theorem \ref{Maintheorem111-2017-irregular}  to end the proof. \qed






{\bf Proof of Theorem \ref{Maintheorem333-2017-alpha-regularsets}. }  If $R_\alpha =X,$ then one can use Theorem \ref{Maintheorem000-2017-banachrecurrence} to complete the proof. Otherwise $I_\alpha\neq \emptyset$ so that
$Int(L_\alpha)\neq \emptyset.$ Denote 
 $\emptyset, QR\cap BR_1 ,BR_1, BR_2,BR_3, BR_4,BR_5$ by $Z_1,Z_2,\cdots,Z_7$ respectively.
 By  Theorem  \ref{Maintheorem222-2017-levelsets} and Proposition \ref{Prop-Fact-sup-H-alpha-a}, for any $i=1,2,\cdots,6,$
$$h_{top}(R_\alpha\cap Z_{i+1}\setminus Z_i)\geq \sup_{a\in Int(L_\alpha)} h_{top}(R_\alpha(a)\cap Z_{i+1}\setminus Z_i)$$$$=\sup_{a\in Int(L_\alpha)} h_{top}(R_\alpha(a) )=\sup_{a\in Int(L_\alpha)}H_\alpha(a)=h_{top}=h_{top}(R_\alpha).\qed$$

\begin{Rem}
By Lemma \ref{Lem-minimalmeasuredense} and  \ref{Lem-Tran-in-BanV}, $Tran\subseteq BV$ so that Theorem \ref{Maintheorem000-2017-banachrecurrence}, Theorem \ref{Maintheorem111-2017-irregular},   Theorem  \ref{Maintheorem222-2017-levelsets} and Theorem \ref{Maintheorem333-2017-alpha-regularsets} can be stated w.r.t. $BV.$

\end{Rem}

\section{Comments and Further Questions}\label{section-commentsquestions}
\subsection{Regularity, Quasiregularity} \label{section-regular----irregular}


  In    \cite{Oxt}  Oxtoby also introduced another several concepts (for quasiregular point, also see \cite{DGS}):
 \begin{eqnarray*}
 QR_{erg}  &:=&\{\text{points generic  for ergodic measures}\}
 = \cup_{\mu\in M_{erg}(T,X)} G_\mu, 
 \\ QR_{d}  &:=&\{\text{points of density in } QR \}= \cup_{\mu\in M(T,X)} (G_\mu\cap S_\mu), 
 \\ R  &:=&\{\text{regular points of } T\}= QR_{d} \cap QR_{erg} = \cup_{\mu\in M_{erg}(T,X)} (G_\mu\cap S_\mu).
 \end{eqnarray*}
Such sets are all $T-$invariant and remark that $$R \subseteq QR_{d} \cup QR_{erg}  \subseteq QR. $$
In \cite{T16} $  R$ and  $  QR_d \setminus   R  =QR_d \setminus   QR_{erg}  $ were considered but  $   QR_{erg} \setminus R=QR_{erg} \setminus QR_d$ and  $  QR \setminus (  QR_{erg}\cup QR_d)$ are not considered.

\begin{Thm}\label{Thm-Oxtoby}
  Suppose that $(X,\,T)$   has  $g$-almost product property.
   Let $\alpha:M(T,X)\rightarrow \mathbb{R}$ be a continuous function  with  $Int(L_\alpha)\neq \emptyset.$ 
 If $(X,\,T)$ is not uniquely ergodic and  there is an invariant measure
with full support,
       then\\
       (1) if $\alpha$ satisfies   $[A.1]$, then   $$h_{top}(BR_1\cap R_\alpha(a)\cap
       QR \setminus (  QR_{erg}\cup QR_d))=h_{top}(R_\alpha(a))\,\,\text{where }\,a\in Int(L_\alpha)\,\,\,\text{ and }$$
        $$ h_{top}(  QR \setminus (  QR_{erg}\cup QR_d))
         =h_{top}(BR_1\cap R_\alpha\cap QR \setminus (  QR_{erg}\cup QR_d))=h_{top}(T).$$
        (2)  $h_{top}(  QR_{erg} \setminus R)
        =h_{top}(BR_1\cap R_\alpha\cap (QR_{erg} \setminus R)) =h_{top}(T).$  
 \end{Thm}
{\bf Proof.} In the proof of  Theorem \ref{mianthm--Maintheorem222-2017-levelsets}, $G_{K_0}$ is in fact contained in  $Tran_1\cap R_\alpha(a)\cap QR \setminus (  QR_{erg}\cup QR_d)$ and $[A.2]$ is not used.
By Lemma \ref{Lem-Tran-in-BR}, $Tran\subseteq BR$ so that $  Tran_1\subseteq BR_1.$ Thus one can follow the proof of Theorem \ref{mianthm--Maintheorem222-2017-levelsets} to get the first result in item (1).

By    Proposition \ref{Main-prop-levelsets-0} and Proposition \ref{Prop-Fact-sup-H-alpha-a}, $h_{top}(BR_1\cap R_\alpha\cap  QR \setminus (  QR_{erg}\cup QR_d)) $
$$ \geq \sup_{a\in Int(L_\alpha)} h_{top}(BR_1\cap R_\alpha(a)\cap QR \setminus (  QR_{erg}\cup QR_d))$$$$=\sup_{a\in Int(L_\alpha)} h_{top}(R_\alpha(a) )=\sup_{a\in Int(L_\alpha)}H_\alpha(a)=h_{top}=h_{top}(R_\alpha). $$

On the other hand, using $  Tran_1\subseteq BR_1$, entropy-dense and transitively-single-saturated, we have  $h_{top}(BR_1\cap R_\alpha\cap QR_{erg} \setminus R)$
$$\geq \sup  \{h_{top}(Tran_1\cap R_\alpha \cap G_\mu \setminus QR_d)|\,\mu\in M_{erg}(T,X),\,S_\mu\neq X\}$$
$$= \sup  \{h_{top}(  G^T_\mu  )|\,\mu\in M_{erg}(T,X),\,S_\mu\neq X\}= \sup  \{h_{\mu}(  T )|\,\mu\in M_{erg}(T,X),\,S_\mu\neq X\}$$$$=h_{top}=h_{top}(R_\alpha).\qed$$

\begin{Thm}\label{Thm-Oxtoby1111111111111}
  Suppose that $(X,\,T)$   has  entropy-dense property.
  Then $h_{top}(  R\setminus Tran)
        = h_{top}(T).$  
 \end{Thm}

 {\bf Proof. } By variational principle and entropy-dense property,  for any $\epsilon>0,$ there is an ergodic
 measure $\mu$ with $S_\mu\neq X$ such that $h_\mu(T)>h_{top}-\epsilon.$ By    \cite[Theorem 3]{Bowen1} $h_{top}(G_\mu\cap S_\mu)=h_\mu(T)$ and note that
 $G_\mu\cap S_\mu\subseteq R\setminus Tran$ so that $h_{top}(  R\setminus Tran)>h_{top}-\epsilon.$ \qed

\subsection{The results on $QW$}

 Let 
 \begin{eqnarray*}   W  &:=&\{\, x\in  QW  \,|\,S_\mu=C_x  \text{ for every }\mu \in M_x \},
   \\   V   &:=&\{ \,x\in  QW  |\,\exists \,\mu\in M_x  \text{ such that } S_\mu=C_x\},\\
   S  &:=&\{ \,x\in  X  |\,\cap_{\mu\in M_x} S_\mu\neq \emptyset\};\\
   QW_1  &:=&  W,\,\,
   QW_2   :=  V\cap S ,\,\,\,\, QW_3   :=  V, \\
   QW_4  &:=& V\cup (QW\cap S) ,\,\,\,\,\, QW_5   :=  QW,
          \end{eqnarray*}


\begin{maintheorem}\label{Maintheorem-recall-adv2016-2017-banachrecurrence}
Suppose that $(X,\,T)$   has  saturated property and  entropy-dense property.
 If $(X,\,T)$ is not uniquely ergodic and  there is an invariant measure
with full support,    then
    $\{\emptyset, QR\cap QW_1 ,QW_1, QW_2,QW_3, QW_4,QW_5 \}$ has   full entropy gap w.r.t. $  Tran $. 
  Similar arguments hold  w.r.t. $Tran\cap I_\alpha,$ $Tran\cap R_\alpha,$ $Tran\cap R_\alpha(a).$

 \end{maintheorem}

 {\bf Proof.} All constructed $K $ in \cite{T16} satisfy $C_K=X$ so that $G_K\subseteq QW\cap Tran$ (in this case $G_K=G^T_K$).   
 Thus  transitively-saturated property is not necessary for this result and saturated property is enough.

For  $\{\emptyset, QR\cap QW_1 ,QW_1, QW_2,QW_3\}$, one can follow the  constructed $K $ in \cite{T16} that  all consist of convex sum of finite measures.
 On the other hand, for $QW_5\setminus QW_4$ and $QW_4\setminus QW_3$, the  constructed $K $ in \cite{T16} is not contained in convex sum of finite measures so that convex-saturated is not enough. Thus here we assume the system to be saturated. \qed

\begin{Rem}
The results of Theorem \ref{Maintheorem000-2017-banachrecurrence}, Theorem \ref{Maintheorem111-2017-irregular},
 Theorem  \ref{Maintheorem222-2017-levelsets},  Theorem \ref{Maintheorem333-2017-alpha-regularsets} and Theorem \ref{Maintheorem-recall-adv2016-2017-banachrecurrence}
  are all restricted on transitive points. By
 their  assumption 
  the system is not minimal so that almost periodic case is not contained in these results. Almost periodic case will be discussed in \cite{DongTian2017}.

\end{Rem}

\subsection{Statistical $\omega-$limit sets} Recently several concepts of statistical $\omega-$limit sets were introduced in \cite{DongTian2016}.
\begin{Def}
 For $x\in X$ and $\xi=\overline{d},   \,  \underline{d},   \,  B^*,  \,   B_*$,   a point $y\in X$ is called $x-\xi-$accessible,   if for any $ \epsilon>0,  \,  N (x,  V_\epsilon (y))\text{ has positive   density w.  r.  t.   }\xi,  $ where     $V_\epsilon (x)$ denotes  the  ball centered at $x$ with radius $\epsilon$.
 Let $$\omega_{\xi}(x):=\{y\in X\,  |\,   y\text{ is } x-\xi-\text{accessible}\}.  $$  For convenience,   it is called {\it $\xi-\omega$-limit set of $x$}.
 \end{Def}






For any $x\in X$, if $\omega_{B_*}(x)=\emptyset,$ then from \cite{DongTian2016} we know that  $x$ satisfies only one of following twelve cases:
  \begin{description}
  \item[Case  (1)    ] \,   $  \emptyset=\omega_{B_*}(x)\subsetneq\omega_{\underline{d}}(x)= \omega_{\overline{d}}(x)= \omega_{B^*}(x)= \omega_f(x);$
  \item[Case  (1')    ]  \,  \,   $  \emptyset=\omega_{B_*}(x)\subsetneq\omega_{\underline{d}}(x)= \omega_{\overline{d}}(x)= \omega_{B^*}(x)\subsetneq\omega_f(x);$

  \item[Case  (2)    ]  \,  \,   $\emptyset=\omega_{B_*}(x)\subsetneq    \omega_{\underline{d}}(x)=
\omega_{\overline{d}}(x)  \subsetneq \omega_{B^*}(x)= \omega_f(x) ;$
\item[Case  (2')    ]  \,  \,   $\emptyset=\omega_{B_*}(x)\subsetneq  \omega_{\underline{d}}(x)=
\omega_{\overline{d}}(x)  \subsetneq \omega_{B^*}(x)\subsetneq \omega_f(x) ;$

  \item[Case  (3)    ]  \,  \,   $    \emptyset=\omega_{B_*}(x)=\omega_{\underline{d}}(x)\subsetneq
\omega_{\overline{d}}(x)= \omega_{B^*}(x)= \omega_f(x) ;$
  \item[Case  (3')    ] \,  \,   $\emptyset=\omega_{B_*}(x)=\omega_{\underline{d}}(x)\subsetneq
\omega_{\overline{d}}(x)= \omega_{B^*}(x)\subsetneq \omega_f(x) ;$
\item[Case  (4)    ]   \,  \,
$\emptyset=\omega_{B_*}(x)\subsetneq  \omega_{\underline{d}}(x)\subsetneq
\omega_{\overline{d}}(x)= \omega_{B^*}(x)= \omega_f(x) ;$
 \item[Case  (4')    ]  \,  \,
$\emptyset=\omega_{B_*}(x)\subsetneq  \omega_{\underline{d}}(x)\subsetneq
\omega_{\overline{d}}(x)= \omega_{B^*}(x)\subsetneq \omega_f(x) ;$

  \item[Case  (5)    ]  \,  \,   $\emptyset=\omega_{B_*}(x)=  \omega_{\underline{d}}(x) \subsetneq \omega_{\overline{d}}(x)\subsetneq \omega_{B^*}(x)= \omega_f(x);$
 \item[Case  (5')    ] \,  \,   $\emptyset=\omega_{B_*}(x)=  \omega_{\underline{d}}(x) \subsetneq \omega_{\overline{d}}(x)\subsetneq \omega_{B^*}(x)\subsetneq \omega_f(x);$
\item[Case  (6)    ]   \,  \,
$\emptyset=\omega_{B_*}(x)\subsetneq \omega_{\underline{d}}(x) \subsetneq \omega_{\overline{d}}(x)\subsetneq \omega_{B^*}(x)= \omega_f(x);$
 \item[Case  (6')    ]  \,  \,
$\emptyset=\omega_{B_*}(x)\subsetneq \omega_{\underline{d}}(x) \subsetneq \omega_{\overline{d}}(x)\subsetneq \omega_{B^*}(x)\subsetneq \omega_f(x).$

  \end{description}

\begin{maintheorem}\label{Maintheorem-recall-adv2016-2017-recurrent-statisticallanguage}
Suppose that $(X,\,T)$   has  transitively-saturated property and  entropy-dense property.
 If $(X,\,T)$ is not uniquely ergodic and  there is an invariant measure
with full support,    then
    $\{x\in Rec|\, x \text{ satisfies Case } i\}$ has   full entropy gap w.r.t. $  Tran,\,i=1,2,\cdots, 6 $. 
Similar arguments holds w.r.t. $Tran\cap I_\alpha,$ $Tran\cap R_\alpha,$ $Tran\cap R_\alpha(a).$

 \end{maintheorem}

{\bf Proof.}
For any $x\in X,   $  from \cite{DongTian2016} we know $\omega_{\underline{d}}(x)= \bigcap_{\mu\in M_x} S_\mu$,
  $\omega_{\overline{d}}(x)=C_x\neq \emptyset,$ $\omega_{B_*}(x)= \bigcap_{\mu\in M^*_x} S_\mu$ and
 $\omega_{B^*}(x)= C_x^*\neq \emptyset.$ Thus $$x\in BR \Leftrightarrow x\in \omega_{B^*}(x)\text{ and } x\in QW \Leftrightarrow x\in \omega_{\overline{d}}(x).$$
  The construction of $x$ in the proof of results for $BR_i$ in present paper  and $QW_i$ in \cite{T16} always satisfies that $x\in Tran\cap BR,$ $\omega_{B^*}(x)= C_x^*=X=\omega_T(x)$ and $M^*_x=M(T,X)$.
    Since the dynamical systems of main theorems  are not minimal but minimal points are dense in the whole space so that for any $x\in Tran$, $\omega_{B_*}(x)=\emptyset.$
     Thus one can check  that $BR_1$ belongs to Case (2), $BR_2\setminus BR_1$ and $BR_4\setminus BR_3$ belong  to Case (6), $BR_3\setminus BR_2$ and $BR_5\setminus BR_4$
      belong  to Case (5), $QW_1$ belongs to Case (1), $QW_2\setminus QW_1$ and $QW_4\setminus QW_3$ belong  to Case (4), $QW_3\setminus QW_2$ and $QW_{5}\setminus QW_4$ belong
       to Case (3). So Cases (1)-(6)
     have full topological entropy $Tran$, $Tran\cap I_\alpha,$ $  Tran\cap R_\alpha,$ $Tran\cap R_\alpha(a) $ respectively by using Theorem \ref{Maintheorem000-2017-banachrecurrence}, Theorem \ref{Maintheorem111-2017-irregular},   Theorem  \ref{Maintheorem222-2017-levelsets},  Theorem \ref{Maintheorem333-2017-alpha-regularsets} and Theorem \ref{Maintheorem-recall-adv2016-2017-banachrecurrence}.\qed

\begin{Rem}
Suppose that $(X,\,T)$   has $g$-product property.
 If $(X,\,T)$ is not uniquely ergodic and  there is an invariant measure
with full support,    then
    $$\{x\in Rec|\, x \text{ satisfies Case } i\}$$ has   full entropy gap w.r.t. $  Tran,\,i=1,2,\cdots, 6 $. 
Similar arguments holds w.r.t. $Tran\cap I_\alpha,$ $Tran\cap R_\alpha,$ $Tran\cap R_\alpha(a).$ Here it is not necessary to assume uniform separation, since   only the constructed $K$ for the gap-sets $QW_4\setminus QW_3$ and  $QW_5\setminus QW_4$  used saturated property, but  transitively-convex-saturated  is enough for others.

 \end{Rem}
\begin{Rem}
Cases (1)-(6) and Cases (1')-(6') restricted on non-recurrent points are considered in \cite{DongTian2016}.
\end{Rem}

\begin{Rem}
The case of $\omega_{B_*}(x)\neq \emptyset$ can also be classified into many cases (including almost periodic case)  restricted on recurrent points or  non-recurrent points which will be studied in \cite{DongTian2017}.
\end{Rem}



 \subsection{Topological entropy on $Rec\setminus BR$}

Let $ Rec^\# :=Rec\setminus BR$,
 \begin{eqnarray*}  Rec_1^\# &:=&\{\, x\in  Rec^\# \,|\, C_x =C^*_x   \}, \,\,
     Rec_2^\#  :=  \{\, x\in  Rec^\# \,|\, C_x \subsetneq C^*_x  \},
   \\ S^\#  &:=&\{ \,x\in  X  |\,\cap_{\mu\in M_x} S_\mu\neq \emptyset\};\\
    W_i^\# &:=&\{\, x\in  Rec_i^\# \,|\,S_\mu=C_x   \text{ for every }\mu \in M_x \},
  \\  V_i^\#  &:=&\{ \,x\in  Rec_i^\#  |\,\exists \,\mu\in M_x   \text{ such that } S_\mu=C_x\},\\
      Rec_{i,1}  &:=&  W_i^\#, \,\,
   Rec_{i,2}   :=  V_i^\#\cap S^\# ,\,\,\,\, Rec_{i,3}   :=  V_i^\#, \\
   Rec_{i,4}  &:=& V_i^\#\cup (Rec^\#\cap S^\#) ,\,\,\,\,\, Rec_{i,5}   :=  Rec_i^\#, \text{ where } i=1,2.
          \end{eqnarray*}

\begin{Que}\label{Question-recurrent-minus-banachrecurrent}
  Suppose that $(X,\,T)$   has  $g$-almost product property and uniform separation.
 If $(X,\,T)$ is not uniquely ergodic and  there is an invariant measure
with full support,  then whether
$\{\emptyset, Rec_{i,1}, Rec_{i,2},Rec_{i,3},Rec_{i,4},Rec_{i,5}\}$ has  full entropy gap w.r.t. $ X$ for any $i=1,2$? Similar questions can be asked  w.r.t. $  I_\alpha,$ $  R_\alpha,$ $  R_\alpha(a).$

 \end{Que}

\begin{Rem}

Under the assumption of Question \ref{Question-recurrent-minus-banachrecurrent}, $Tran\subseteq BR$ so that $Rec_{i,j}\cap Tran=\emptyset.$

\end{Rem}

It is also an open question that the  sets of $Rec_{i,j}$ and their gap-sets are nonempty in what kind of dynamical systems.
  For $Rec_{1,1}$, it is nonempty in full shifts over finite symbols by following theorem. Note that any full shift  over finite symbols satisfies the assumption of Question \ref{Question-recurrent-minus-banachrecurrent} so that Question \ref{Question-recurrent-minus-banachrecurrent} is possibly  meaningful.

  \begin{Thm} Let $\sigma:\{0,1,\cdots,k-1\}^\mathbb{N}\rightarrow \{0,1,\cdots,k-1\}^\mathbb{N}$ be a full shift where $k\geq 2.$ Then $Rec_{1,1}\neq \emptyset.$

 \end{Thm}

{\bf Proof.} We only need to consider $k=2.$
We learned a example from \cite{HZ}  that  there is  a topological mixing but uniquely ergodic subshift $\Lambda$ of two symbols for which the unique invariant measure  is supported on a fixed point so that in this example $\emptyset \neq Tran(\sigma|_\Lambda)\subseteq Rec_{1,1}\cap  \{x|\,\emptyset\neq \omega_{B_*}(x)=\omega_{\underline{d}}(x)= \omega_{\overline{d}}(x)= \omega_{B^*}(x)\subsetneq \omega_f(x)\}.$ \qed 

\medskip

 The  sets of $Rec_{i,j}\cap \{x|\, \omega_{B_*}(x)= \emptyset\}$ are related to Cases (1')-(6') restricted on recurrent points.  So similar question can be asked:

\begin{Que}\label{Question-recurrent-minus-banachrecurrent---2}
  Suppose that $(X,\,T)$   has  $g$-almost product property and uniform separation.
 If $(X,\,T)$ is not uniquely ergodic and  there is an invariant measure
with full support,  then whether Cases (1')-(6') restricted on recurrent points all carry full topological entropy?
 \end{Que}

\section*{Acknowlegements}  The research of Y. Huang was supported by National Natural Science Foundation of China (Grant No.11771459).  The research of X. Tian was  supported by National Natural Science Foundation of China  (grant no. 11671093). The research of X. Wang was supported by National Natural Science Foundation of China, Tian Yuan Special Foundation(Grant No.11426198) and the Natural Science Foundation of Guangdong Province,China (Grant No.2015A030310166).


\end{document}